\documentclass[12pt]{amsart}

\usepackage{amssymb}
\usepackage{amsmath}
\usepackage{amsthm}
\usepackage{amsfonts}

\usepackage{a4wide}

\newtheorem{theorem}{Theorem}
\newtheorem{corollary}{Corollary}

\theoremstyle{definition}
\newtheorem{definition}{Definition}

\newtheorem{conjecture}{Conjecture}[section]

\numberwithin{equation}{section}

\title[Harmonically balanced capitulation]{Cyclic cubic number fields \\ with harmonically balanced capitulation}

\author{Bill Allombert}
\address{Universit\'e de Bordeaux\\351 cours de la Lib\'eration\\33405 Bordeaux\\France}
\email{Bill.Allombert@math.u-bordeaux.fr}

\author{Daniel C. Mayer}
\address{Naglergasse 53\\8010 Graz\\Austria}
\email{algebraic.number.theory@algebra.at}
\urladdr{http://www.algebra.at}

\thanks{Research of second author supported by the Austrian Science Fund (FWF): projects J0497-PHY, P26008-N25, and by the Research Executive Agency of the European Union (EUREA)}



\date{Wednesday 26 July 2023}

\begin{document}


\begin{abstract}
It is proved that
\(c=689\,347=31\cdot 37\cdot 601\)
is the smallest conductor of
a cyclic cubic number field \(K\)
whose maximal unramified pro-\(3\)-extension
\(E=\mathrm{F}_3^\infty(K)\)
possesses an automorphism group
\(G=\mathrm{Gal}(E/K)\)
of order \(6561\)
with coinciding relation and generator rank
\(d_2(G)=d_1(G)=3\)
and harmonically balanced transfer kernels
\(\varkappa(G)\in S_{13}\).
\end{abstract}

\maketitle


\section{Introduction}
\label{s:Intro}

\noindent
Let \(p\) be a prime number.
Finite \(p\)-groups \(G\) with \textit{balanced presentation},
that is, with relation rank \(d_2(G)\) equal to the generator rank \(d_1(G)\),
\begin{equation}
\label{eqn:Schur}
d_2(G)=
\mathrm{dim}_{\mathbb{F}_p}H^2(G,\mathbb{F}_p)=
\mathrm{dim}_{\mathbb{F}_p}H^1(G,\mathbb{F}_p)=
d_1(G),
\end{equation}
have attracted the vigilance and interest of researchers since the beginning.
Issai Schur dubbed such groups \textit{closed},
but today it is more usual to call them \textit{Schur groups}.

In \(1964\),
one year before his famous joint paper on infinite class field towers with Golod,
Shafarevich
\cite{Sh1964}
established an application of closed groups
of the greatest importance and with fundamental impact on class field theory.
Let \(K\) be an algebraic number field
with \(p\)-class rank \(\varrho=\mathrm{rank}_p\mathrm{Cl}(K)\),
signature \((r_1,r_2)\),
and torsion free Dirichlet unit rank \(r=r_1+r_2-1\).
Let \(\theta\in\lbrace 0,1\rbrace\) be an indicator for the existence of
a primitive \(p\)-th root of unity
\(\zeta=\exp(2\pi\sqrt{-1}/p)\) in \(K\).
Then the relation rank \(d_2(G)\) of the Galois group
\(G=\mathrm{Gal}(\mathrm{F}_p^\infty(K)/K)\)
of the maximal unramified pro-\(p\)-extension \(\mathrm{F}_p^\infty(K)\) of \(K\)
is bounded
\begin{equation}
\label{eqn:Shafarevich}
\varrho\le d_2(G)\le\varrho+r+\theta,
\end{equation}
according to Shafarevich
\cite[Thm. 5.1, p. 28]{Ma2015c}.
In particular,
imaginary quadratic fields \(k=\mathbb{Q}(\sqrt{d})\)
with negative discriminant \(d<0\),
which have the simplest possible signature \((0,1)\)
(except \((1,0)\) for the rational number field \(\mathbb{Q}\))
require a \textit{Schur \(\sigma\)-group}
\(\mathrm{Gal}(\mathrm{F}_p^\infty(k)/k)\)
for all \(p\)-class towers with
an odd prime \(p\ge 3\)
\cite{Ag1998}.
Shafarevich himself immediately drew some conclusions
\cite[pp. 91--92]{Sh1964}
about the few non-abelian \(3\)-class towers
for \(d=-4027\) and \(d=-3299\)
which were known at this early stage,
due to Scholz and Taussky
\cite{SoTa1934}.
Unfortunately, his interpretations of these \(3\)-class towers
for \(3\)-class groups with abelian type invariants \((3,3)\) and \((9,3)\) 
were incorrect
(his claimed groups were only Schur but not Schur \(\sigma\)),
and his important theorem on the bounds for \(d_2(G)\) remained largely unnoticed
for about thirty years.

In \(1975\), several events happened precipitately.
On the one hand, Andozhskii and Tsvetkov
\cite{AnTs1974},
\cite{An1975}
found the first closed finite \(3\)-groups \(G\)
with elementary tricyclic commutator quotient \(G/G^\prime\simeq (3,3,3)\).
But on the other hand,
Koch and Venkov
\cite{KoVe1975}
proved that, for odd primes \(p\ge 3\),
an imaginary quadratic number field
\(k=\mathbb{Q}(\sqrt{d})\) with
\(p\)-class group \(\mathrm{Cl}_p(k)\)
of \(p\)-rank \(\varrho\ge 3\) possesses 
an unbounded \(p\)-class field tower.
In particular, for \(\mathrm{Cl}_3(k)\simeq (3,3,3)\)
the Galois group \(\mathrm{Gal}(F_3^\infty(k)/k)\)
of the \(3\)-class tower 
must be an infinite Schur \(\sigma\)-group
and cannot be one of the Andozhskii-Tsvetkov groups
(briefly \textit{AT-groups}).

This is exactly the point where our present article sets in.
In \S\
\ref{s:Andozhskii},
we investigate AT-groups more closely,
determine their order \(3^8=6561\)
and position in the descendant tree
\cite{Ma2016},
and discover with surprise that all of them possess
\textit{harmonically balanced capitulation} (HBC),
that is, their transfer kernels
\(\varkappa(G)=(\ker(T_j))_{j=1}^{13}\)
are cyclic of order \(3\)
and form a permutation in the symmetric group \(S_{13}\)
on thirteen letters.

Then we try to realize AT-groups as automorphism groups
of maximal unramified pro-\(3\)-extensions of
suitable algebraic number fields \(K\),
necessarily different from imaginary quadratic fields,
by the result of Koch and Venkov.
Since real quadratic fields of type \((3,3,3)\)
are firstly very sparse
and secondly also require a \(\sigma\)-group
with generator- and relator-inverting automorphism
(though with looser bounds \(d_1(G)\le d_2(G)\le d_1(G)+1\)),
they disqualify for the realization of AT-groups.
In \S\
\ref{s:CyclicCubic},
we employ the next simplest number fields
with signature \((3,0)\) and unit rank \(r=2\),
but still absolutely Galois over \(\mathbb{Q}\),
that is, the cyclic cubic fields \(K\),
which quite frequently have
elementary tricyclic \(3\)-class groups
\(\mathrm{Cl}_3(K)\simeq (3,3,3)\),
provided their conductor \(c\) is divisible by three primes,
whence they arise as quartets \((K_1,K_2,K_3,K_4)\)
sharing a common discriminant \(d=c^2\),
as known from earlier works by G. Gras \(1973\)
\cite{Gr1973}
and Ayadi \(1995\)
\cite{Ay1995},
\cite{AAI2001}.

Indeed, we find a rather sparse sequence of conductors \(c=q_1q_2q_3\),
exclusively with graph 2 of category I,
\(q_2\leftarrow q_1\rightarrow q_3\),
in the sense of G. Gras and Ayadi,
for which precisely one component \(K:=K_1\)
has \(\mathrm{Cl}_3(K)\simeq (3,3,3)\) and HBC,
whereas the other three components \(K_2,K_3,K_4\) of the quartet
have \(\mathrm{Cl}_3(K_i)\simeq (3,3)\) and
\(\varkappa(K_i)=(1243)\), type G.16.
It is published as sequence A359310
in the On-line Encyclopedia of Integer Sequences (OEIS).

A drawback of cyclic cubic fields \(K\) is the broad spread
for the relation rank \(d_2(G)\) of the \(3\)-class tower group
\(G=\mathrm{Gal}(\mathrm{F}_3^\infty(K)/K)\),
\begin{equation}
\label{eqn:D2}
3=\varrho=d_1(G)\le d_2(G)\le\varrho+r+\theta=3+2+0=5,
\end{equation}
which enables the occurrence of numerous smaller groups \(G\)
of order \(3^6=729\) as parents of AT-groups with \(d_2(G)=5\),
or of order \(3^7=2187\) as siblings of AT-groups with \(d_2(G)=4\),
instead of the desired minimal relation rank \(d_2(G)=3\).
We had to wait for the \(25\)-th term \(c=\mathbf{689\,347}\)
of the OEIS sequence A359310
\cite{OEIS2023}
until a proper AT-group occurred,
and this delay required lots of CPU time and RAM storage
on most powerful super computers with multiple processor threads,
unimaginable even on bigger workstations.


\section{Closed Andozhskii-Tsvetkov groups}
\label{s:Andozhskii}

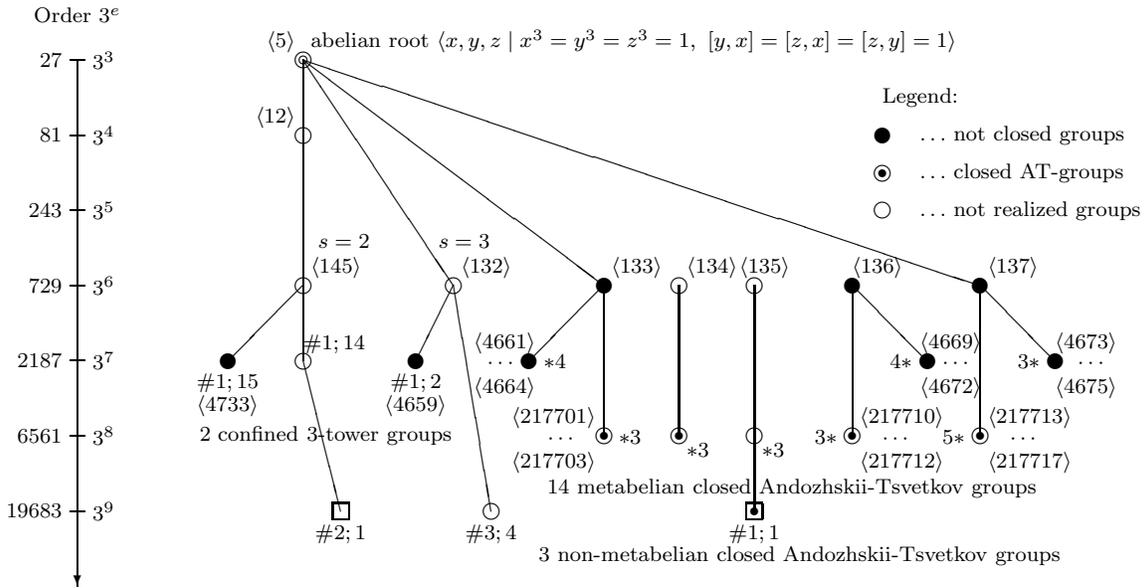
\begin{figure}[ht]
\caption{Tree of \(3\)-groups \(G\) with \(G/G^\prime\simeq (3,3,3)\)}
\label{fig:Groups333}

{\tiny

\setlength{\unitlength}{1.0cm}
\begin{picture}(10,8)(-7.5,-7)


\put(-9,0.5){\makebox(0,0)[cb]{Order \(3^e\)}}
\put(-9,0){\line(0,-1){7}}
\multiput(-9.1,0)(0,-1){7}{\line(1,0){0.2}}
\put(-9.2,0){\makebox(0,0)[rc]{\(27\)}}
\put(-8.8,0){\makebox(0,0)[lc]{\(3^3\)}}
\put(-9.2,-1){\makebox(0,0)[rc]{\(81\)}}
\put(-8.8,-1){\makebox(0,0)[lc]{\(3^4\)}}
\put(-9.2,-2){\makebox(0,0)[rc]{\(243\)}}
\put(-8.8,-2){\makebox(0,0)[lc]{\(3^5\)}}
\put(-9.2,-3){\makebox(0,0)[rc]{\(729\)}}
\put(-8.8,-3){\makebox(0,0)[lc]{\(3^6\)}}
\put(-9.2,-4){\makebox(0,0)[rc]{\(2187\)}}
\put(-8.8,-4){\makebox(0,0)[lc]{\(3^7\)}}
\put(-9.2,-5){\makebox(0,0)[rc]{\(6561\)}}
\put(-8.8,-5){\makebox(0,0)[lc]{\(3^8\)}}
\put(-9.2,-6){\makebox(0,0)[rc]{\(19683\)}}
\put(-8.8,-6){\makebox(0,0)[lc]{\(3^9\)}}
\put(-9,-6){\vector(0,-1){1}}


\put(1.7,-0.5){\makebox(0,0)[lc]{Legend:}}
\put(1.7,-1.0){\circle*{0.2}}
\put(2.2,-1.0){\makebox(0,0)[lc]{\(\ldots\) not closed groups}}
\put(1.7,-1.5){\circle{0.2}}
\put(1.7,-1.5){\circle*{0.1}}
\put(2.2,-1.5){\makebox(0,0)[lc]{\(\ldots\) closed AT-groups}}
\put(1.7,-2.0){\circle{0.2}}
\put(2.2,-2.0){\makebox(0,0)[lc]{\(\ldots\) not realized groups}}

\put(-6.1,0.1){\makebox(0,0)[rb]{\(\langle 5\rangle\)}}
\put(-5.9,0.1){\makebox(0,0)[lb]{abelian root \(\langle x,y,z\mid
x^3=y^3=z^3=1,\ \lbrack y,x\rbrack=\lbrack z,x\rbrack=\lbrack z,y\rbrack=1\rangle\)}}
\put(-6,0){\circle{0.2}}
\put(-6,0){\circle{0.1}}

\put(-6,0){\line(0,-1){1}}
\put(-6.1,-0.9){\makebox(0,0)[rb]{\(\langle 12\rangle\)}}
\put(-6,-1){\circle{0.2}}
\put(-6,-1){\line(0,-1){2}}
\put(-5.8,-2.5){\makebox(0,0)[lb]{\(s=2\)}}
\put(-5.9,-2.9){\makebox(0,0)[lb]{\(\langle 145\rangle\)}}
\put(-6,-3){\circle{0.2}}
\put(-6,-3){\line(-1,-1){1}}
\put(-7,-4){\circle*{0.2}}
\put(-7,-4.3){\makebox(0,0)[cc]{\(\#1;15\)}}
\put(-7,-4.6){\makebox(0,0)[cc]{\(\langle 4733\rangle\)}}
\put(-6,-3){\line(0,-1){1}}
\put(-6,-3.9){\makebox(0,0)[lb]{\(\#1;14\)}}
\put(-6,-4){\circle{0.2}}
\put(-6,-4){\line(1,-4){0.5}}
\put(-5.6,-6.1){\framebox(0.2,0.2){}}
\put(-5.5,-6.3){\makebox(0,0)[cc]{\(\#2;1\)}}

\put(-5.7,-5){\makebox(0,0)[cc]{\(2\) confined \(3\)-tower groups}}

\put(-6,0){\line(2,-3){2}}
\put(-4.2,-2.5){\makebox(0,0)[lb]{\(s=3\)}}
\put(-3.9,-2.9){\makebox(0,0)[lb]{\(\langle 132\rangle\)}}
\put(-4,-3){\circle{0.2}}
\put(-4,-3){\line(-1,-2){0.5}}
\put(-4.5,-4){\circle*{0.2}}
\put(-4.5,-4.3){\makebox(0,0)[cc]{\(\#1;2\)}}
\put(-4.5,-4.6){\makebox(0,0)[cc]{\(\langle 4659\rangle\)}}
\put(-4,-3){\line(1,-6){0.5}}
\put(-3.5,-6){\circle{0.2}}
\put(-3.5,-6.3){\makebox(0,0)[cc]{\(\#3;4\)}}

\put(-6,0){\line(4,-3){4}}
\put(-1.9,-2.9){\makebox(0,0)[lb]{\(\langle 133\rangle\)}}
\put(-2,-3){\circle*{0.2}}
\put(-2,-3){\line(0,-1){2}}
\put(-2,-5){\circle{0.2}}
\put(-2,-5){\circle*{0.1}}
\put(-1.8,-5.1){\makebox(0,0)[lb]{\(\ast 3\)}}
\put(-2.1,-4.9){\makebox(0,0)[rb]{\(\langle 217701\rangle\)}}
\put(-2.4,-5.1){\makebox(0,0)[rb]{\(\cdots\)}}
\put(-2.1,-5.5){\makebox(0,0)[rb]{\(\langle 217703\rangle\)}}

\put(-2,-3){\line(-1,-1){1}}
\put(-3,-4){\circle*{0.2}}
\put(-2.8,-4.1){\makebox(0,0)[lb]{\(\ast 4\)}}
\put(-2.9,-3.9){\makebox(0,0)[rb]{\(\langle 4661\rangle\)}}
\put(-3.2,-4.1){\makebox(0,0)[rb]{\(\cdots\)}}
\put(-2.9,-4.5){\makebox(0,0)[rb]{\(\langle 4664\rangle\)}}

\put(-0.9,-2.9){\makebox(0,0)[lb]{\(\langle 134\rangle\)}}
\put(-1,-3){\circle{0.2}}
\put(-1,-3){\line(0,-1){2}}
\put(-1,-5){\circle{0.2}}
\put(-1,-5){\circle*{0.1}}
\put(-0.9,-5.1){\makebox(0,0)[lt]{\(\ast 3\)}}

\put(-0.2,-2.9){\makebox(0,0)[lb]{\(\langle 135\rangle\)}}
\put(0,-3){\circle{0.2}}
\put(0,-3){\line(0,-1){2}}
\put(0.5,-5.7){\makebox(0,0)[cc]{\(14\) metabelian closed Andozhskii-Tsvetkov groups}}
\put(0,-5){\circle{0.2}}
\put(0.1,-5.1){\makebox(0,0)[lt]{\(\ast 3\)}}
\put(0,-5){\line(0,-1){1}}
\put(-0.1,-6.1){\framebox(0.2,0.2){}}
\put(0,-6){\circle*{0.1}}
\put(0,-6.3){\makebox(0,0)[cc]{\(\#1;1\)}}
\put(0.6,-6.6){\makebox(0,0)[cc]{\(3\) non-metabelian closed Andozhskii-Tsvetkov groups}}

\put(1.3,-2.9){\makebox(0,0)[lb]{\(\langle 136\rangle\)}}
\put(1.3,-3){\circle*{0.2}}
\put(1.3,-3){\line(0,-1){2}}
\put(1.3,-5){\circle{0.2}}
\put(1.3,-5){\circle*{0.1}}
\put(1.1,-5.1){\makebox(0,0)[rb]{\(3\ast\)}}
\put(1.4,-4.9){\makebox(0,0)[lb]{\(\langle 217710\rangle\)}}
\put(1.7,-5.1){\makebox(0,0)[lb]{\(\cdots\)}}
\put(1.4,-5.5){\makebox(0,0)[lb]{\(\langle 217712\rangle\)}}

\put(1.3,-3){\line(1,-1){1}}
\put(2.3,-4){\circle*{0.2}}
\put(2.1,-4.1){\makebox(0,0)[rb]{\(4\ast\)}}
\put(2.2,-3.9){\makebox(0,0)[lb]{\(\langle 4669\rangle\)}}
\put(2.5,-4.1){\makebox(0,0)[lb]{\(\cdots\)}}
\put(2.2,-4.5){\makebox(0,0)[lb]{\(\langle 4672\rangle\)}}

\put(-6,0){\line(3,-1){9}}
\put(3.1,-2.9){\makebox(0,0)[lb]{\(\langle 137\rangle\)}}
\put(3,-3){\circle*{0.2}}
\put(3,-3){\line(0,-1){2}}
\put(3,-5){\circle{0.2}}
\put(3,-5){\circle*{0.1}}
\put(2.8,-5.1){\makebox(0,0)[rb]{\(5\ast\)}}
\put(3.1,-4.9){\makebox(0,0)[lb]{\(\langle 217713\rangle\)}}
\put(3.4,-5.1){\makebox(0,0)[lb]{\(\cdots\)}}
\put(3.1,-5.5){\makebox(0,0)[lb]{\(\langle 217717\rangle\)}}

\put(3,-3){\line(1,-1){1}}
\put(4,-4){\circle*{0.2}}
\put(3.8,-4.1){\makebox(0,0)[rb]{\(3\ast\)}}
\put(4,-3.9){\makebox(0,0)[lb]{\(\langle 4673\rangle\)}}
\put(4.3,-4.1){\makebox(0,0)[lb]{\(\cdots\)}}
\put(4,-4.5){\makebox(0,0)[lb]{\(\langle 4675\rangle\)}}


\end{picture}
}
\end{figure}

\noindent
According to Koch and Venkov
\cite{KoVe1975},
\textit{Schur \(\sigma\)-groups} \(S\) are known to be
mandatory for realizations
\(S\simeq\mathrm{Gal}(\mathrm{F}_p^\infty(k)/k)\)
by \(p\)-class field towers of \textit{imaginary} quadratic fields \(k\),
with an odd prime \(p\).
They possess a balanced presentation \(d_1(S)=d_2(S)\) with coinciding
generator rank \(d_1(S)=\dim_{\mathbb{F}_p}H^1(S,\mathbb{F}_p)\)
and relation rank \(d_2(S)=\dim_{\mathbb{F}_p}H^2(S,\mathbb{F}_p)\),
and an automorphism \(\sigma\in\mathrm{Aut}(S)\)
acting as inversion \(x\mapsto x^{-1}\)
on the commutator quotient \(S/\lbrack S,S\rbrack\).
However, in the older literature,
for instance Shafarevich
\cite[\S\ 6, pp. 88--91]{Sh1964},
there also appear \textit{Schur groups} with balanced presentation,
but without a generator inverting \(\sigma\)-automorphism,
and they are called \textit{closed},
according to the original terminology by Schur.
In the present article,
we are interested in finite closed \(3\)-groups
given by Andozhskii and Tsvetkov (\textit{AT-groups})
\cite{AnTs1974,An1975},
whose position in the descendant tree of 
\(3\)-groups \(G\) with elementary tricyclic
commutator quotient \(G/\lbrack G,G\rbrack\simeq (3,3,3)\)
is illuminated in Figure
\ref{fig:Groups333}.


In \S\
\ref{ss:AndozhskiiIdentified},
we identify  the \(17\) \textit{closed Andozhskii-Tsvetkov groups}
as the smallest \(3\)-groups of type \((3,3,3)\) with balanced presentation.
Their order is either \(3^8=6561\) or \(3^9=19683\).
In \S\
\ref{ss:RealizationAndozhskii},
we show that three or four of them can be realized as Galois groups
of the \(3\)-class tower of \textit{cyclic cubic fields}.
Incidentally,
we point out that related but not closed \(3\)-groups
(see \(j=2\) in Table
\ref{tbl:HBC})
are realized by numerous \(3\)-class field towers over 
\textit{totally complex \(S_3\)-fields} \(K\),
which are unramified extensions of imaginary quadratic fields
\(k=\mathbb{Q}(\sqrt{d})\)
with \(3\)-class group \(\mathrm{Cl}_3(k)\simeq (3,3)\),
capitulation type \(\mathrm{H}.4\), \(\varkappa(k)\sim (4111)\),
and three abelian type invariants of rank \(3\) in \(\alpha(k)\sim (111,111,111,21)\).
The latter \(3\)-class towers are \textit{confined}.


\subsection{Identification of closed Andozhskii-Tsvetkov groups}
\label{ss:AndozhskiiIdentified}

\noindent
First we prove \textit{existence} and determine their \textit{number} (Thm.
\ref{thm:Andozhskii}).
Then we compute \textit{invariants} (Cor.
\ref{cor:Andozhskii}
and
\ref{cor:Taussky}).

\begin{theorem}
\label{thm:Andozhskii}
Among the finite \(3\)-groups \(G\)
with commutator quotient \(G/G^\prime\simeq (3,3,3)\),
there exist precisely \(\mathbf{14}\) \textbf{metabelian closed}
groups \(S\) of order \(\#(S)=3^8\)
with identifiers
\begin{equation}
\label{eqn:Metab}
S\simeq\langle 6561,217700+i\rangle \text{ where } 1\le i\le 6 \text{ or } 10\le i\le 17,
\end{equation}
and \(\mathbf{3}\) \textbf{non-metabelian closed}
groups \(S\) of order \(\#(S)=3^9\)
with identifiers
\begin{equation}
\label{eqn:NonMetab}
S\simeq\langle 6561,217700+i\rangle-\#1;1 \text{ where } 7\le i\le 9.
\end{equation}
They possess a trivial Schur multiplier \(M(S)=H_2(S,\mathbb{Q}/\mathbb{Z})=0\)
and a balanced presentation \(d_1(S)=d_2(S)\)
with coinciding generator rank
\(d_1(S)=\dim_{\mathbb{F}_p}H^1(S,\mathbb{F}_p)\)
and relation rank
\(d_2(S)=\dim_{\mathbb{F}_p}H^2(S,\mathbb{F}_p)\).
The class is \(\mathrm{cl}(S)=3\) for soluble length \(\mathrm{sl}(S)=2\)
and \(\mathrm{cl}(S)=4\) for \(\mathrm{sl}(S)=3\).
They possess harmonically balanced transfer kernels \(\varkappa(S)\in S_{13}\).
\end{theorem}

\begin{proof}
By a search in the SmallGroups database
\cite{BEO2005},
extended to order \(3^9\) by the \(p\)-group generation algorithm
\cite{Nm1977,Ob1990,HEO2005},
the finite closed Andozhskii-Tsvetkov \(3\)-groups \(S\)
are identified.
There are no hits of order \(\#(S)\le 3^7\),
\(14\) hits of order \(\#(S)=3^8\),
and only three hits of order \(\#(S)=3^9\).
The non-metabelian groups are characterized by
their relative identifiers defined in the ANUPQ package
\cite{GNO2006}.
See Figure
\ref{fig:Groups333}
and Table
\ref{tbl:HBC}.
\end{proof}


\begin{corollary}
\label{cor:Andozhskii}
Each of the \(17\) closed Andozhskii groups \(S\) in Theorem
\ref{thm:Andozhskii}
shares a common \textbf{Artin pattern} \((\varkappa,\alpha)\)
with its ancestor \(A\simeq\langle 729,130+j\rangle\),
as given in Table
\ref{tbl:HBC},
where
\begin{equation}
\label{eqn:Connection}
j=
\begin{cases}
3 & \text{ for } 1\le i \le 3, \\
4 & \text{ for } 4\le i \le 6, \\
5 & \text{ for } 7\le i \le 9, \\
6 & \text{ for } 10\le i \le 12, \\
7 & \text{ for } 13\le i \le 17. \\
\end{cases}
\end{equation}
\end{corollary}

\begin{proof}
According to the \textit{theorem on the antitony} of the Artin pattern
\cite[\S\S\ 5.1--5.4, pp. 78--87]{Ma2016},
it suffices to calculate the \textit{stable} transfer kernels
of the five ancestors \(A\) of the \(17\) closed groups in Theorem
\ref{thm:Andozhskii}.
They are of order \(\#(A)=3^6\)
and have much simpler presentations.
It turns out that the transfer kernel type (TKT) is \textit{harmonically balanced},
that is, a permutation in the symmetric group \(S_{13}\)
of degree \(13\).
\end{proof}


\noindent
We recall that the transfer kernel (capitulation kernel) \(\ker(T_i)\)
of an Artin transfer homomorphism \(T_i:\,G/G^\prime\to M_i/M_i^\prime\)
\cite{Ma2016}
from a \(3\)-group \(G\) with \(G/G^\prime\simeq (3,3,3)\)
to one of its \(13\) maximal subgroups \(M_i\), \(1\le i\le 13\),
is called of \textit{Taussky type} A, if the meet
\(\ker(T_i)\cap M_i>1\) is non-trivial,
and of \textit{Taussky type} B,
if \(\ker(T_i)\cap M_i=1\) is trivial.

\begin{corollary}
\label{cor:Taussky}
The Artin patterns \((\varkappa(A),\alpha(A))\) of the six groups
\(A=\langle 729,130+j\rangle\), \(2\le j\le 7\),
share the common property that the Taussky type of the transfer kernels
\(\varkappa(A)_i=\ker(T_i)\)
is determined uniquely by the \textbf{AQI} \(M_i/M_i^\prime\) of the corresponding maximal subgroup \(M_i\):
\begin{equation}
\label{eqn:Taussky}
\begin{aligned}
\alpha(A)_i=M_i/M_i^\prime\simeq (211) & \Longleftrightarrow \varkappa(A)_i\cap M_i>1, \text{ Taussky type }\mathrm{A}, \\
\alpha(A)_i=M_i/M_i^\prime\simeq (22)  & \Longleftrightarrow \varkappa(A)_i\cap M_i=1, \text{ Taussky type }\mathrm{B},
\end{aligned}
\end{equation}
for all \(1\le i\le 13\). Here the Abelian Quotient Invariants are written in logarithmic form.
\end{corollary}

\begin{proof}
This law follows by comparing the \(1\)-dimensional transfer kernels (lines) \(\ker(T_i)\) in Table
\ref{tbl:HBC}
to the \(2\)-dimensional subspaces (planes) \(M_j\) of \(S/S^\prime\simeq (3,3,3)\),
\(1\le i,j\le 13\).
\end{proof}


\renewcommand{\arraystretch}{1.1}

\begin{table}[ht]
\caption{TKT \(\varkappa\), AQI \(\alpha\), ranks \(\rho\), and operator of ancestors \(\langle 729,130+j\rangle\)}
\label{tbl:HBC}
\begin{center}
{\scriptsize
\begin{tabular}{|c||c|c|c|c|c|c|c|c|c|c|c|c|c||c||c|}
\hline
 \(j\) & \multicolumn{13}{|c||}{\(\varkappa\) respectively \(\alpha\)} & \(\rho\) & Action \\
\hline
 \(2\) &  1 &  2 &  3 &   6 & 11 &  9 & 10 &   4 & 13 &   5 & 12 &  8 &   7 &  & \(\langle 24,12\rangle\) \\
       & 22 & 22 & 22 & 211 & 22 & 22 & 22 & 211 & 22 & 211 & 22 & 22 & 211 & \(3^4,2^9\) & \\
\hline
 \(3\) &  9 &  2 &  3 &   6 & 10 &   8 &   4 &  11 & 12 &  13 &   5 &  1 &   7 &  & \(\langle 6,2\rangle\) \\
       & 22 & 22 & 22 & 211 & 22 & 211 & 211 & 211 & 22 & 211 & 211 & 22 & 211 & \(3^7,2^6\) & \\
\hline
 \(4\) &  1 &  7 &  3 &  2 & 10 &  9 &   4 &  11 &   5 &  12 & 13 &  8 &  6 &  & \(\langle 4,1\rangle\) \\
       & 22 & 22 & 22 & 22 & 22 & 22 & 211 & 211 & 211 & 211 & 22 & 22 & 22 & \(3^4,2^9\) & \\
\hline
 \(5\) &  9 &  7 &  3 &  2 &   4 &   8 & 11 &  10 & 13 &   5 & 12 &  1 &  6 &  & \(\langle 3,1\rangle\) \\
       & 22 & 22 & 22 & 22 & 211 & 211 & 22 & 211 & 22 & 211 & 22 & 22 & 22 & \(3^4,2^9\)  & \\
\hline
 \(6\) & 12 &  7 &  3 &  2 &  9 &  5 &  8 &  1 & 10 & 11 &  4 &  13 &  6 &  & \(\langle 6,2\rangle\) \\
       & 22 & 22 & 22 & 22 & 22 & 22 & 22 & 22 & 22 & 22 & 22 & 211 & 22 & \(3^1,2^{12}\) & \\
\hline
 \(7\) & 10 &  7 &  3 &   6 &   8 &  4 &  1 &  9 & 13 &   5 & 12 &  11 &  2 &  & \(\langle 24,3\rangle\) \\
       & 22 & 22 & 22 & 211 & 211 & 22 & 22 & 22 & 22 & 211 & 22 & 211 & 22 & \(3^4,2^9\) & \\
\hline
\end{tabular}
}
\end{center}
\end{table}


\subsection{Second order invariants of candidate groups with HBC}
\label{ss:Candidates}

\noindent
In Table
\ref{tbl:AQI2HBC},
we present all possible candidates for metabelian \(3\)-groups \(G\)
with HBC, characterized uniquely by their absolute identifier
\(\langle 3^{\mathrm{lo}},\mathrm{id}\rangle\) in the SmallGroups database
\cite{BEO2005},
\cite{MAGMA6561}
with order \(3^{\mathrm{lo}}\), logarithmic order \(6\le\mathrm{lo}\le 8\)
and numerical identifier \(\mathrm{id}\).
Crucial invariants of these groups are the nuclear rank \(\nu\),
the \(p\)-multiplicator rank \(\mu=d_2\), which coincides with the \textit{relation rank},
and the number of descendants \(N_s\)
and of capable descendants \(C_s\) for all possible step sizes \(1\le s\le\nu\).
The most important invariant, however,
indispensable for the unambiguous identification
(\(\varkappa\) and \(\alpha\) are insufficient),
and demanding extreme computational challenge for the number theoretic verification,
is the \textit{Artin pattern of second order}
\cite{Ma2020},
\begin{equation}
\label{eqn:AQI2GT}
\alpha^{(2)}:=\lbrack H_j/H_j^\prime;(H_{j,\ell}/H_{j,\ell}^\prime)_{\ell=1}^{n_j}\rbrack_{j=1}^{13},
\end{equation}
consisting of \textit{logarithmic abelian quotient invariants} (AQI)
of all maximal subgroups \(H_j\), \(1\le j\le 13\),
and second maximal subgroups \(H_{j,\ell}\), \(1\le j\le 13\), with
\(1\le\ell\le 4\), \(n_j=4\), when \(H_j/H_j^\prime\simeq (22)\hat{=}(9,9)\), and
\(1\le\ell\le 13\), \(n_j=13\), when \(H_j/H_j^\prime\simeq (211)\hat{=}(9,3,3)\).

The group \(\langle 729,132\rangle\) is forbidden
for cyclic cubic fields,
because of its relation rank \(6\).
The group \(\langle 729,134\rangle\) and all its descendants are forbidden
for cyclic cubic fields,
because they have a wrong action by \(\langle 4,1\rangle\),
which is only allowed for \textit{cyclic quartic fields}.
Among the groups of order \(3^6=729\),
two with identifiers \(\mathrm{id}\in\lbrace 135,137\rbrace\)
cannot be distinguished by second order invariants \(\alpha^{(2)}\).
Similarly for groups of order \(3^7=2187\),
two with identifiers \(\mathrm{id}\in\lbrace 4662,4663\rbrace\) and
two with identifiers \(\mathrm{id}\in\lbrace 4673,4675\rbrace\)
cannot be separated by \(\alpha^{(2)}\).
Similarly for groups of order \(3^8=6561\),
two with identifiers \(\mathrm{id}\in\lbrace 217702,217703\rbrace\) and
two with identifiers \(\mathrm{id}\in\lbrace 217713,217717\rbrace\)
cannot be separated by \(\alpha^{(2)}\).

So far,
three or four AT-groups of order \(3^8=6561\) have been \textbf{realized} as
\(3\)-class tower groups of cyclic cubic fields.
Either \(\langle 6561,217702\rangle\) or \(\langle 6561,217703\rangle\)
by the conductor \(c=\mathbf{1\,406\,551}\),
either \(\langle 6561,217713\rangle\) or \(\langle 6561,217717\rangle\)
by the conductors \(c=\mathbf{689\,347}\) and \(c=\mathbf{869\,611}\), and
\(\langle 6561,217710\rangle\) unambiguously
by \(c=\mathbf{753\,787}\) and \(c=\mathbf{796\,779}\).


\renewcommand{\arraystretch}{1.1}

\begin{table}[ht]
\caption{Second order invariants and propagation of \(3\)-groups with HBC}
\label{tbl:AQI2HBC}
\begin{center}
{\tiny
\begin{tabular}{|c|r|l|c|c|c|}
\hline
 lo & id     & \(\alpha^{(2)}\) & \(\nu\) & \(\mu\) & \((N_s/C_s)_{s=1}^\nu\) \\
\hline
  6 & \(\mathbf{136}\) & \(\lbrack (21^2);(21^2)^4(21)^9\rbrack^1\lbrack (2^2);(21^2)^4\rbrack^{12}\) & 2 & 5 & (4/0;3/0) \\
\hline
  6 &    132 & \(\lbrack (21^2);(21^2)^4(21)^9\rbrack^4\lbrack (2^2);(21^2)^4\rbrack^9\)    & 3 & 6 & (3/0;6/1;4/4) \\   
  6 &    134 & \(\lbrack (21^2);(21^2)^4(21)^9\rbrack^4\lbrack (2^2);(21^2)^4\rbrack^9\)    & 2 & 5 & (2/0;3/0) \\
  6 &    135 & \(\lbrack (21^2);(21^2)^4(21)^9\rbrack^4\lbrack (2^2);(21^2)^4\rbrack^9\)    & 2 & 5 & (2/0;3/3) \\
  6 & \(\mathbf{137}\) & \(\lbrack (21^2);(21^2)^4(21)^9\rbrack^4\lbrack (2^2);(21^2)^4\rbrack^9\)    & 2 & 5 & (3/0;5/0) \\
\hline
  6 & \(\mathbf{133}\) & \(\lbrack (21^2);(21^2)^4(21)^9\rbrack^7\lbrack (2^2);(21^2)^4\rbrack^6\)    & 2 & 5 & (4/0;3/0) \\
\hline
  7 &   4669 & \(\lbrack (21^2);(21^2)^4(21)^9\rbrack^1\lbrack (2^2);(31^2)(21^2)^3\rbrack^4\lbrack (2^2);(21^2)^4\rbrack^8\) & 0 & 4 & --- \\
  7 & \(\mathbf{4670}\) & \(\lbrack (21^2);(2^21)(21^2)^{12}\rbrack^1\lbrack (2^2);(2^21)(21^2)^3\rbrack^3\lbrack (2^2);(21^2)^4\rbrack^9\) & 0 & 4 & --- \\
  7 &   4671 & \(\lbrack (21^2);(2^21)(21^2)^3(2^2)^9\rbrack^1\lbrack (2^2);(2^21)(21^2)^3\rbrack^3\lbrack (2^2);(21^2)^4\rbrack^9\) & 0 & 4 & --- \\
  7 &   4672 & \(\lbrack (21^2);(2^21)(21^2)^3(2^2)^9\rbrack^1\lbrack (2^2);(2^21)(21^2)^3\rbrack^3\lbrack (2^2);(21^2)^4\rbrack^9\) & 0 & 4 & --- \\
\hline
  7 & \(\mathbf{4673}\) & \(\lbrack (21^2);(21^2)^4(21)^9\rbrack^3\lbrack (21^2);(2^21)(21^2)^3(2^2)^9\rbrack^1\lbrack (2^2);(2^21)(21^2)^3\rbrack^3\lbrack (2^2);(21^2)^4\rbrack^6\) & 0 & 4 & --- \\
  7 &   4674 & \(\lbrack (21^2);(21^2)^4(21)^9\rbrack^3\lbrack (21^2);(2^21)(21^2)^{12}\rbrack^1\lbrack (2^2);(2^21)(21^2)^3\rbrack^3\lbrack (2^2);(21^2)^4\rbrack^6\) & 0 & 4 & --- \\
  7 & \(\mathbf{4675}\) & \(\lbrack (21^2);(21^2)^4(21)^9\rbrack^3\lbrack (21^2);(2^21)(21^2)^3(2^2)^9\rbrack^1\lbrack (2^2);(2^21)(21^2)^3\rbrack^3\lbrack (2^2);(21^2)^4\rbrack^6\) & 0 & 4 & --- \\
\hline
  7 &   4661 & \(\lbrack (21^2);(21^2)^4(21)^9\rbrack^6\lbrack (21^2);(2^21)(21^2)^{12}\rbrack^1\lbrack (2^2);(2^21)(21^2)^3\rbrack^3\lbrack (2^2);(21^2)^4\rbrack^3\) & 0 & 4 & --- \\
  7 & \(\mathbf{4662}\) & \(\lbrack (21^2);(21^2)^4(21)^9\rbrack^6\lbrack (21^2);(2^21)(21^2)^3(2^2)^9\rbrack^1\lbrack (2^2);(2^21)(21^2)^3\rbrack^3\lbrack (2^2);(21^2)^4\rbrack^3\) & 0 & 4 & --- \\
  7 & \(\mathbf{4663}\) & \(\lbrack (21^2);(21^2)^4(21)^9\rbrack^6\lbrack (21^2);(2^21)(21^2)^3(2^2)^9\rbrack^1\lbrack (2^2);(2^21)(21^2)^3\rbrack^3\lbrack (2^2);(21^2)^4\rbrack^3\) & 0 & 4 & --- \\
  7 &   4664 & \(\lbrack (21^2);(21^2)^4(21)^9\rbrack^5\lbrack (21^2);(31^2)(21^2)^3(31)^9\rbrack^2\lbrack (2^2);(31^2)(21^2)^3\rbrack^2\lbrack (2^2);(21^2)^4\rbrack^4\) & 0 & 4 & --- \\
\hline
  8 & \(\mathbf{217710}\) & \(\lbrack (21^2);(2^21)(21^2)^{12}\rbrack^1\lbrack (2^2);(31^2)^3(2^21)\rbrack^1\lbrack (2^2);(2^21)(21^2)^3\rbrack^2\lbrack (2^2);(31^2)(21^2)^3\rbrack^9\) & 0 & 3 & --- \\
  8 & 217711 & \(\lbrack (21^2);(2^21)(21^2)^3(2^2)^9\rbrack^1\lbrack (2^2);(31^2)^3(2^21)\rbrack^1\lbrack (2^2);(2^21)(21^2)^3\rbrack^2\lbrack (2^2);(31^2)(21^2)^3\rbrack^9\) & 0 & 3 & --- \\
  8 & 217712 & \(\lbrack (21^2);(2^21)(21^2)^3(2^2)^9\rbrack^1\lbrack (2^2);(31^2)^3(2^21)\rbrack^1\lbrack (2^2);(2^21)(21^2)^3\rbrack^2\lbrack (2^2);(31^2)(21^2)^3\rbrack^9\) & 0 & 3 & --- \\
\hline
  8 & \(\mathbf{217713}\) & \(\lbrack (21^2);(2^21)(21^2)^3(2^2)^9\rbrack^4\lbrack (2^2);(2^21)(21^2)^3\rbrack^8\lbrack (2^2);(21^2)^4\rbrack^1\) & 0 & 3 & --- \\
  8 & 217714 & \(\lbrack (21^2);(2^21)(21^2)^3(2^2)^9\rbrack^2\lbrack (21^2);(2^21)(21^2)^{12}\rbrack^2\lbrack (2^2);(2^21)(21^2)^3\rbrack^8\lbrack (2^2);(21^2)^4\rbrack^1\) & 0 & 3 & --- \\
  8 & 217715 & \(\lbrack (21^2);(2^21)(21^2)^3(2^2)^9\rbrack^3\lbrack (21^2);(2^21)(21^2)^{12}\rbrack^1\lbrack (2^2);(2^21)(21^2)^3\rbrack^8\lbrack (2^2);(21^2)^4\rbrack^1\) & 0 & 3 & --- \\
  8 & 217716 & \(\lbrack (21^2);(2^21)(21^2)^3(2^2)^9\rbrack^1\lbrack (21^2);(2^21)(21^2)^{12}\rbrack^3\lbrack (2^2);(2^21)(21^2)^3\rbrack^8\lbrack (2^2);(21^2)^4\rbrack^1\) & 0 & 3 & --- \\
  8 & \(\mathbf{217717}\) & \(\lbrack (21^2);(2^21)(21^2)^3(2^2)^9\rbrack^4\lbrack (2^2);(2^21)(21^2)^3\rbrack^8\lbrack (2^2);(21^2)^4\rbrack^1\) & 0 & 3 & --- \\
\hline
  8 & 217701 & \(\lbrack (21^2);(31^2)(21^2)^3(31)^9\rbrack^6\lbrack (21^2);(2^21)(21^2)^{12}\rbrack^1\) & 0 & 3 & --- \\
    &        & \(\lbrack (2^2);(2^21)(21^2)^3\rbrack^2\lbrack (2^2);(31^2)(21^2)^3\rbrack^3\lbrack (2^2);(31^2)(2^21)^3\rbrack^1\) & & & \\
  8 & \(\mathbf{217702}\) & \(\lbrack (21^2);(31^2)(21^2)^3(31)^9\rbrack^6\lbrack (21^2);(2^21)(21^2)^3(2^2)^9\rbrack^1\) & 0 & 3 & --- \\
    &        & \(\lbrack (2^2);(2^21)(21^2)^3\rbrack^2\lbrack (2^2);(31^2)(21^2)^3\rbrack^3\lbrack (2^2);(31^2)(2^21)^3\rbrack^1\) & & & \\
  8 & \(\mathbf{217703}\) & \(\lbrack (21^2);(31^2)(21^2)^3(31)^9\rbrack^6\lbrack (21^2);(2^21)(21^2)^3(2^2)^9\rbrack^1\) & 0 & 3 & --- \\
    &        & \(\lbrack (2^2);(2^21)(21^2)^3\rbrack^2\lbrack (2^2);(31^2)(21^2)^3\rbrack^3\lbrack (2^2);(31^2)(2^21)^3\rbrack^1\) & & & \\
\hline

\end{tabular}
}
\end{center}
\end{table}


\section{Realization of groups with HBC by algebraic number fields}
\label{s:CyclicCubic}


\subsection{Realization as 3-class field tower groups}
\label{ss:RealizationAndozhskii}

\noindent
Since the groups in Theorem
\ref{thm:Andozhskii}
are non-\(\sigma\) groups,
they cannot be realized by any quadratic field,
neither imaginary nor real. Therefore,
we investigated the possible Galois actions (Table
\ref{tbl:HBC})
on the five ancestors \(A=\mathrm{SmallGroup}(729,130+j)\).
It turned out that the unique non-metabelian case
\(j=5\) can only be realized by \textit{cyclic cubic} fields,
\(j=4\) by cyclic quartic fields, and
\(j\in\lbrace 3,6,7\rbrace\) by \textit{cyclic cubic or sextic} fields.
We show that certain metabelian descendants \(S\) for \(j\in\lbrace 3,6,7\rbrace\)
can actually be realized as Galois groups
\(\mathrm{Gal}(\mathrm{F}_3^\infty(K)/K)\simeq S\)
of maximal unramified pro-\(3\)-extensions
of cyclic cubic fields \(K\) with \(53\) conductors
\(c\) in the OEIS sequence A359310
\cite{OEIS2023}, 
from \(\mathbf{59\,031}\) to \(\mathbf{1\,406\,551}\),
and with \(3\)-class group \(\mathrm{Cl}_3(K)\simeq (3,3,3)\).


\begin{theorem}
\label{thm:Main1}
If a number field \(K/\mathbb{Q}\)
with elementary tricyclic \(3\)-class group \(\mathrm{Cl}_3(K)\simeq (3,3,3)\)
possesses the Artin pattern \((\varkappa(K),\alpha(K))\)
with harmonically balanced capitulation type 
\(\varkappa(K)\sim (9,2,3,6,10,8,4,11,12,13,5,1,7)\)
and  abelian type invariants
\(\alpha(K)\sim ((22)^{3},211,22,(211)^3,22,(211)^2,22,211)\),
then \(K/\mathbb{Q}\) must be cyclic cubic or sextic,
and has a metabelian \(3\)-class field tower with automorphism group
\(\mathrm{Gal}(\mathrm{F}_3^\infty(K)/K)\simeq\)
\begin{equation}
\label{eqn:Main1}
\langle 3^8,217700+i\rangle,\ 1\le i\le 3,
\text{ or } \langle 3^7,4660+k\rangle,\ 1\le k\le 4,
\text{ or } \langle 3^6,133\rangle.
\end{equation}
\end{theorem}


\begin{theorem}
\label{thm:Main2}
If a number field \(K/\mathbb{Q}\)
with elementary tricyclic \(3\)-class group \(\mathrm{Cl}_3(K)\simeq (3,3,3)\)
possesses the Artin pattern \((\varkappa(K),\alpha(K))\)
with harmonically balanced capitulation type 
\(\varkappa(K)\sim (12,7,3,2,9,5,8,1,10,11,4,13,6)\)
and abelian type invariants
\(\alpha(K)\sim ((22)^{11},211,22)\),
then \(K/\mathbb{Q}\) must be cyclic cubic or sextic,
and has a metabelian \(3\)-class field tower with automorphism group
\(\mathrm{Gal}(\mathrm{F}_3^\infty(K)/K)\simeq\)
\begin{equation}
\label{eqn:Main2}
\langle 3^8,217700+i\rangle,\ 10\le i\le 12,
\text{ or } \langle 3^7,4669+k\rangle,\ 0\le k\le 3,
\text{ or } \langle 3^6,136\rangle.
\end{equation}
\end{theorem}

\begin{proof}
Theorems
\ref{thm:Main1}
and
\ref{thm:Main2}
are immediate consequences of the Tables
\ref{tbl:HBC}
and
\ref{tbl:AQI2HBC}.
\end{proof}


\begin{conjecture}
\label{cnj:Main3}
If a number field \(K/\mathbb{Q}\)
with elementary tricyclic \(3\)-class group \(\mathrm{Cl}_3(K)\simeq (3,3,3)\)
possesses the Artin pattern \((\varkappa(K),\alpha(K))\)
with harmonically balanced capitulation type 
\(\varkappa(K)\sim (10,7,3,6,8,4,1,9,13,5,12,11,2)\)
and abelian type invariants
\(\alpha(K)\sim ((22)^3,(211)^2,(22)^4,211,22,211,22)\),
then \(K/\mathbb{Q}\) must be cyclic cubic or sextic,
and has a metabelian \(3\)-class field tower with automorphism group
\(\mathrm{Gal}(\mathrm{F}_3^\infty(K)/K)\simeq\)
\begin{equation}
\label{eqn:Main3}
\langle 3^8,217700+i\rangle,\ 13\le i\le 15,
\text{ or } \langle 3^7,4673+k\rangle,\ 0\le k\le 2,
\text{ or } \langle 3^6,137\rangle.
\end{equation}
\end{conjecture}

\noindent
For this situation,
the Tables
\ref{tbl:HBC}
and
\ref{tbl:AQI2HBC}
admit descendants of
\(\langle 3^6,132\rangle\) and \(\langle 3^6,135\rangle\)
as additional candidates.
But, so far, experience provides evidence
that no such realizations occur.
Therefore we conjecture that
this tendency will continue.


\subsection{Invariants of realizing cyclic cubic fields}
\label{ss:Invariants}

\noindent
In Table
\ref{tbl:ConductorsHBC},
we present all \(37\) conductors \(c<10^6\), and a few \(c>10^6\), of the OEIS sequence A359310
\cite{OEIS2023},
and their prime factors (admitting also the prime power \(3^2=9\)).
They give rise to quartets of cyclic cubic number fields \((K_1,K_2,K_3,K_4)\)
with \(3\)-class groups \(\mathrm{Cl}_3(K_1)\simeq (3,3,3)\) and
\(\mathrm{Cl}_3(K_i)\simeq (3,3)\) for \(2\le i\le 4\).
The \textit{rank distribution} of the first order Artin pattern of \(K:=K_1\)
with respect to the \(13\) unramified cyclic cubic extensions \(E_j/K\) is given by
\begin{equation}
\label{eqn:Rho}
\rho=(\ \mathrm{rank}_3(\mathrm{Cl}_3(E_j))\ )_{j=1}^{13},
\end{equation}
denoted with symbolic exponents which indicate iteration.
Finally, the unique or ambiguous candidate for the
metabelian \(3\)-class tower group
\(S=\mathrm{Gal}(\mathrm{F}_3^\infty(K)/K)\)
is given by its absolute identifier in the SmallGroups database
\cite{BEO2002},
\cite{BEO2005},
\cite{MAGMA6561}
(a vertical bar \(\vert\) means \lq\lq or\rq\rq).


In April 2002,
we used the \textit{Voronoi algorithm}
\cite{Vo1896}
and the \textit{Euler product} method
in order to compute the \(15851\) cyclic cubic fields \(K\)
with conductors \(c_{K/\mathbb{Q}}<10^5\)
and their class numbers \(1\le h_K\le 1953\).
Among the fields,
\(4785\) occur as singlets,
\(7726\) in doublets,
\(3132\) in quartets,
and \(208\) in octets.
Twenty years later, in July 2022,
we have confirmed these results,
extended by the class group structures \(\mathrm{Cl}(K)\).
The cyclic cubic fields \(K\) were constructed 
as \textit{ray class fields} over the rational number field,
using Fieker's class field theoretic routines
\cite{Fi2001}
in MAGMA
\cite{MAGMA2023}.
Additionally, we constructed the \(13\)
\textit{unramified cyclic cubic relative extensions} \(E_j/K\)
of absolute degree \(9\),
whenever the \(3\)-class group of \(K\) was \(\mathrm{Cl}_3(K)\simeq (3,3,3)\),
which was the primary goal for the reconstruction
\cite{Ma2022}
in view of the intended realization of Andozhskii-Tsvetkov (AT-)groups with HBC.


\renewcommand{\arraystretch}{1.1}

\begin{table}[ht]
\caption{Conductors of cyclic cubic number fields with HBC}
\label{tbl:ConductorsHBC}
\begin{center}
{\tiny
\begin{tabular}{|r||r|l||c|c|c|}
\hline
 No. & \(c\)           & Factors         & \(\rho\)         & \(\mathrm{Gal}(\mathrm{F}_3^\infty(K)/K)\) & Thm./Cnj. \\
\hline
   1 &     \(59\,031\) & \(3^2,7,937\)   & \((3^1,2^{12})\) & \(\langle 2187,4670\rangle\)      & \ref{thm:Main2} \\
   2 &    \(209\,853\) & \(3^2,7,3331\)  & \((3^7,2^6)\)    & \(\langle 729,133\rangle\)        & \ref{thm:Main1} \\
   3 &    \(247\,437\) & \(3^2,19,1447\) & \((3^1,2^{12})\) & \(\langle 729,136\rangle\)        & \ref{thm:Main2} \\
   4 &    \(263\,017\) & \(19,109,127\)  & \((3^7,2^6)\)    & \(\langle 729,133\rangle\)        & \ref{thm:Main1} \\
   5 &    \(271\,737\) & \(3^2,109,277\) & \((3^4,2^9)\)    & \(\langle 729,137\rangle\)        & \ref{cnj:Main3} \\
   6 &    \(329\,841\) & \(3^2,67,547\)  & \((3^1,2^{12})\) & \(\langle 729,136\rangle\)        & \ref{thm:Main2} \\
   7 &    \(377\,923\) & \(7,13,4153\)   & \((3^1,2^{12})\) & \(\langle 729,136\rangle\)        & \ref{thm:Main2} \\
   8 &    \(407\,851\) & \(37,73,151\)   & \((3^7,2^6)\)    & \(\langle 729,133\rangle\)        & \ref{thm:Main1} \\
   9 &    \(412\,909\) & \(7,61,967\)    & \((3^1,2^{12})\) & \(\langle 729,136\rangle\)        & \ref{thm:Main2} \\
  10 &    \(415\,597\) & \(7,13,4567\)   & \((3^1,2^{12})\) & \(\langle 2187,4670\rangle\)      & \ref{thm:Main2} \\
  11 &    \(416\,241\) & \(3^2,7,6607\)  & \((3^1,2^{12})\) & \(\langle 729,136\rangle\)        & \ref{thm:Main2} \\
  12 &    \(416\,727\) & \(3^2,19,2437\) & \((3^4,2^9)\)    & \(\langle 2187,4673|4675\rangle\) & \ref{cnj:Main3} \\
  13 &    \(462\,573\) & \(3^2,103,499\) & \((3^7,2^6)\)    & \(\langle 2187,4662|4663\rangle\) & \ref{thm:Main1} \\
  14 &    \(474\,561\) & \(3^2,67,787\)  & \((3^4,2^9)\)    & \(\langle 729,137\rangle\)        & \ref{cnj:Main3} \\
  15 &    \(487\,921\) & \(7,43,1621\)   & \((3^4,2^9)\)    & \(\langle 2187,4673|4675\rangle\) & \ref{cnj:Main3} \\
  16 &    \(493\,839\) & \(3^2,37,1483\) & \((3^1,2^{12})\) & \(\langle 2187,4670\rangle\)      & \ref{thm:Main2} \\
  17 &    \(547\,353\) & \(3^2,61,997\)  & \((3^4,2^9)\)    & \(\langle 2187,4673|4675\rangle\) & \ref{cnj:Main3} \\
  18 &    \(586\,963\) & \(13,163,277\)  & \((3^1,2^{12})\) & \(\langle 729,136\rangle\)        & \ref{thm:Main2} \\
  19 &    \(612\,747\) & \(3^2,103,661\) & \((3^7,2^6)\)    & \(\langle 729,133\rangle\)        & \ref{thm:Main1} \\
  20 &    \(613\,711\) & \(7,73,1201\)   & \((3^4,2^9)\)    & \(\langle 729,137\rangle\)        & \ref{cnj:Main3} \\
  21 &    \(615\,663\) & \(3^2,67,1021\) & \((3^4,2^9)\)    & \(\langle 729,137\rangle\)        & \ref{cnj:Main3} \\
  22 &    \(622\,063\) & \(13,109,439\)  & \((3^1,2^{12})\) & \(\langle 2187,4670\rangle\)      & \ref{thm:Main2} \\ 
  23 &    \(648\,427\) & \(13,31,1609\)  & \((3^4,2^9)\)    & \(\langle 729,137\rangle\)        & \ref{cnj:Main3} \\
  24 &    \(651\,829\) & \(37,79,223\)   & \((3^4,2^9)\)    & \(\langle 729,137\rangle\)        & \ref{cnj:Main3} \\
  25 & \(\mathbf{689\,347}\)    & \(31,37,601\)   & \((3^4,2^9)\)    & \(\langle \mathbf{6561},217713|217717\rangle\) & \ref{cnj:Main3} \\
  26 &    \(690\,631\) & \(19,163,223\)  & \((3^4,2^9)\)    & \(\langle 729,137\rangle\)        & \ref{cnj:Main3} \\
  27 & \(\mathbf{753\,787}\)    & \(19,97,409\)   & \((3^1,2^{12})\) & \(\langle \mathbf{6561},217710\rangle\)        & \ref{thm:Main2} \\
  28 & \(\mathbf{796\,779}\)    & \(3^2,223,397\) & \((3^1,2^{12})\) & \(\langle \mathbf{6561},217710\rangle\)        & \ref{thm:Main2} \\
  29 &    \(811\,069\) & \(7,109,1063\)  & \((3^1,2^{12})\) & \(\langle 729,136\rangle\)        & \ref{thm:Main2} \\
  30 &    \(818\,217\) & \(3^2,229,397\) & \((3^1,2^{12})\) & \(\langle 729,136\rangle\)        & \ref{thm:Main2} \\
  31 & \(\mathbf{869\,611}\)    & \(19,37,1237\)  & \((3^4,2^9)\)    & \(\langle \mathbf{6561},217713|217717\rangle\) & \ref{cnj:Main3} \\
  32 &    \(914\,263\) & \(7,211,619\)   & \((3^7,2^6)\)    & \(\langle 729,133\rangle\)        & \ref{thm:Main1} \\
  33 &    \(915\,439\) & \(7,19,6883\)   & \((3^4,2^9)\)    & \(\langle 2187,4673|4675\rangle\) & \ref{cnj:Main3} \\
  34 &    \(922\,167\) & \(3^2,79,1297\) & \((3^1,2^{12})\) & \(\langle 729,136\rangle\)        & \ref{thm:Main2} \\
  35 &    \(936\,747\) & \(3^2,7,14869\) & \((3^7,2^6)\)    & \(\langle 2187,4662|4663\rangle\) & \ref{thm:Main1} \\
  36 &    \(977\,409\) & \(3^2,223,487\) & \((3^4,2^9)\)    & \(\langle 729,137\rangle\)        & \ref{cnj:Main3} \\
  37 &    \(997\,087\) & \(7,13,10957\)  & \((3^4,2^9)\)    & \(\langle 729,137\rangle\)        & \ref{cnj:Main3} \\
\hline
  40 & \(1\,083\,607\) & \(7,283,547\)   & \((3^7,2^6)\)    & \(\langle 729,133\rangle\)        & \ref{thm:Main1} \\
  43 & \(1\,181\,971\) & \(7,19,8887\)   & \((3^7,2^6)\)    & \(\langle 729,133\rangle\)        & \ref{thm:Main1} \\
  50 & \(1\,323\,007\) & \(7,331,571\)   & \((3^7,2^6)\)    & \(\langle 729,133\rangle\)        & \ref{thm:Main1} \\
  53 & \(\mathbf{1\,406\,551}\) & \(19,181,409\)  & \((3^7,2^6)\)    & \(\langle \mathbf{6561},217702|217703\rangle\) & \ref{thm:Main1} \\
\hline
\end{tabular}
}
\end{center}
\end{table}

\newpage

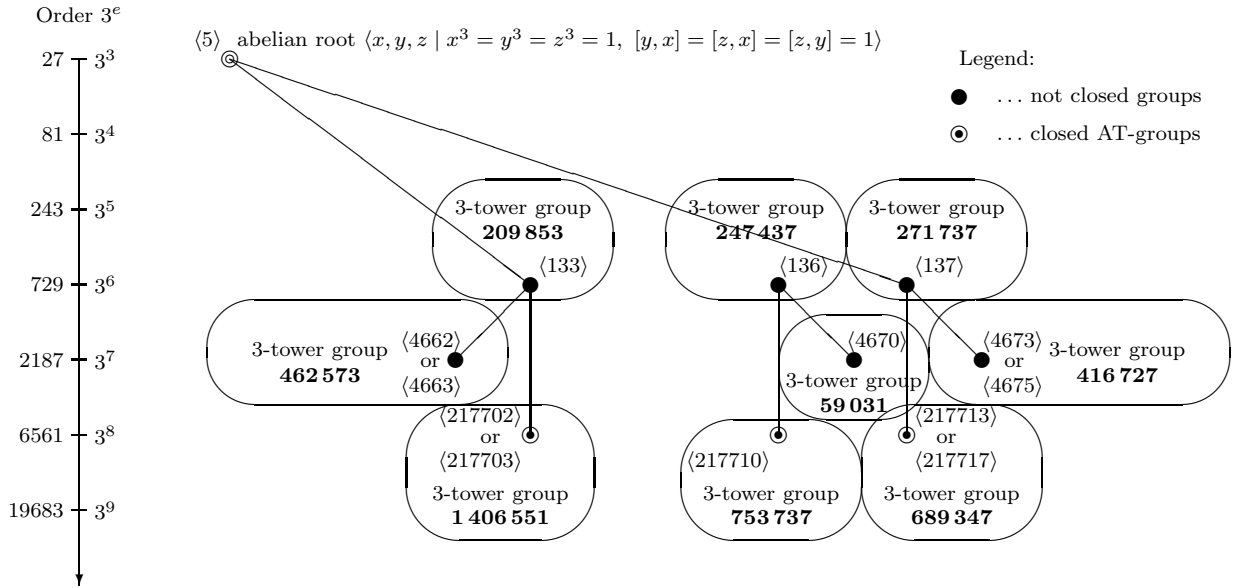
\begin{figure}[ht]
\caption{Realization of groups \(G\) with HBC as \(\mathrm{Gal}(\mathrm{F}_3^\infty(K)/K)\)}
\label{fig:Realizations}

{\tiny

\setlength{\unitlength}{1.0cm}
\begin{picture}(10,8)(-5.9,-7)


\put(-8,0.5){\makebox(0,0)[cb]{Order \(3^e\)}}
\put(-8,0){\line(0,-1){7}}
\multiput(-8.1,0)(0,-1){7}{\line(1,0){0.2}}
\put(-8.2,0){\makebox(0,0)[rc]{\(27\)}}
\put(-7.8,0){\makebox(0,0)[lc]{\(3^3\)}}
\put(-8.2,-1){\makebox(0,0)[rc]{\(81\)}}
\put(-7.8,-1){\makebox(0,0)[lc]{\(3^4\)}}
\put(-8.2,-2){\makebox(0,0)[rc]{\(243\)}}
\put(-7.8,-2){\makebox(0,0)[lc]{\(3^5\)}}
\put(-8.2,-3){\makebox(0,0)[rc]{\(729\)}}
\put(-7.8,-3){\makebox(0,0)[lc]{\(3^6\)}}
\put(-8.2,-4){\makebox(0,0)[rc]{\(2187\)}}
\put(-7.8,-4){\makebox(0,0)[lc]{\(3^7\)}}
\put(-8.2,-5){\makebox(0,0)[rc]{\(6561\)}}
\put(-7.8,-5){\makebox(0,0)[lc]{\(3^8\)}}
\put(-8.2,-6){\makebox(0,0)[rc]{\(19683\)}}
\put(-7.8,-6){\makebox(0,0)[lc]{\(3^9\)}}
\put(-8,-6){\vector(0,-1){1}}


\put(3.7,0){\makebox(0,0)[lc]{Legend:}}
\put(3.7,-0.5){\circle*{0.2}}
\put(4.2,-0.5){\makebox(0,0)[lc]{\(\ldots\) not closed groups}}
\put(3.7,-1){\circle{0.2}}
\put(3.7,-1){\circle*{0.1}}
\put(4.2,-1){\makebox(0,0)[lc]{\(\ldots\) closed AT-groups}}

\put(-6.1,0.1){\makebox(0,0)[rb]{\(\langle 5\rangle\)}}
\put(-5.9,0.1){\makebox(0,0)[lb]{abelian root \(\langle x,y,z\mid
x^3=y^3=z^3=1,\ \lbrack y,x\rbrack=\lbrack z,x\rbrack=\lbrack z,y\rbrack=1\rangle\)}}
\put(-6,0){\circle{0.2}}
\put(-6,0){\circle{0.1}}

\put(-6,0){\line(4,-3){4}}
\put(-1.9,-2.9){\makebox(0,0)[lb]{\(\langle 133\rangle\)}}
\put(-2,-3){\circle*{0.2}}
\put(-2,-3){\line(0,-1){2}}
\put(-2,-5){\circle{0.2}}
\put(-2,-5){\circle*{0.1}}
\put(-2.1,-4.9){\makebox(0,0)[rb]{\(\langle 217702\rangle\)}}
\put(-2.4,-5.1){\makebox(0,0)[rb]{or}}
\put(-2.1,-5.5){\makebox(0,0)[rb]{\(\langle 217703\rangle\)}}

\put(-2.4,-5.5){\oval(2.5,1.8)}
\put(-2.4,-5.8){\makebox(0,0)[cc]{\(3\)-tower group}}
\put(-2.4,-6.1){\makebox(0,0)[cc]{\(\mathbf{1\,406\,551}\)}}

\put(-2,-3){\line(-1,-1){1}}
\put(-3,-4){\circle*{0.2}}
\put(-2.9,-3.9){\makebox(0,0)[rb]{\(\langle 4662\rangle\)}}
\put(-3.2,-4.1){\makebox(0,0)[rb]{or}}
\put(-2.9,-4.5){\makebox(0,0)[rb]{\(\langle 4663\rangle\)}}

\put(-4.3,-3.9){\oval(4,1.4)}
\put(-4.8,-4.0){\makebox(0,0)[cb]{\(3\)-tower group}}
\put(-4.8,-4.3){\makebox(0,0)[cb]{\(\mathbf{462\,573}\)}}

\put(-2.1,-2.4){\oval(2.4,1.6)}
\put(-2.1,-2.0){\makebox(0,0)[cc]{\(3\)-tower group}}
\put(-2.1,-2.3){\makebox(0,0)[cc]{\(\mathbf{209\,853}\)}}

\put(1.0,-2.4){\oval(2.4,1.6)}
\put(1.0,-2.0){\makebox(0,0)[cc]{\(3\)-tower group}}
\put(1.0,-2.3){\makebox(0,0)[cc]{\(\mathbf{247\,437}\)}}

\put(3.4,-2.4){\oval(2.4,1.6)}
\put(3.4,-2.0){\makebox(0,0)[cc]{\(3\)-tower group}}
\put(3.4,-2.3){\makebox(0,0)[cc]{\(\mathbf{271\,737}\)}}

\put(2.3,-4.1){\oval(2,1.4)}
\put(2.3,-4.3){\makebox(0,0)[cc]{\(3\)-tower group}}
\put(2.3,-4.6){\makebox(0,0)[cc]{\(\mathbf{59\,031}\)}}

\put(1.3,-2.9){\makebox(0,0)[lb]{\(\langle 136\rangle\)}}
\put(1.3,-3){\circle*{0.2}}
\put(1.3,-3){\line(0,-1){2}}
\put(1.3,-5){\circle{0.2}}
\put(1.3,-5){\circle*{0.1}}
\put(1.2,-5.5){\makebox(0,0)[rb]{\(\langle 217710\rangle\)}}

\put(1.2,-5.6){\oval(2.4,1.6)}
\put(1.2,-5.8){\makebox(0,0)[cc]{\(3\)-tower group}}
\put(1.2,-6.1){\makebox(0,0)[cc]{\(\mathbf{753\,737}\)}}

\put(1.3,-3){\line(1,-1){1}}
\put(2.3,-4){\circle*{0.2}}
\put(2.2,-3.9){\makebox(0,0)[lb]{\(\langle 4670\rangle\)}}

\put(-6,0){\line(3,-1){9}}
\put(3.1,-2.9){\makebox(0,0)[lb]{\(\langle 137\rangle\)}}
\put(3,-3){\circle*{0.2}}
\put(3,-3){\line(0,-1){2}}
\put(3,-5){\circle{0.2}}
\put(3,-5){\circle*{0.1}}
\put(3.1,-4.9){\makebox(0,0)[lb]{\(\langle 217713\rangle\)}}
\put(3.4,-5.1){\makebox(0,0)[lb]{or}}
\put(3.1,-5.5){\makebox(0,0)[lb]{\(\langle 217717\rangle\)}}

\put(3.6,-5.5){\oval(2.4,1.8)}
\put(3.6,-5.8){\makebox(0,0)[cc]{\(3\)-tower group}}
\put(3.6,-6.1){\makebox(0,0)[cc]{\(\mathbf{689\,347}\)}}

\put(5.3,-3.9){\oval(4,1.4)}
\put(5.8,-4.0){\makebox(0,0)[cb]{\(3\)-tower group}}
\put(5.8,-4.3){\makebox(0,0)[cb]{\(\mathbf{416\,727}\)}}

\put(3,-3){\line(1,-1){1}}
\put(4,-4){\circle*{0.2}}
\put(4,-3.9){\makebox(0,0)[lb]{\(\langle 4673\rangle\)}}
\put(4.3,-4.1){\makebox(0,0)[lb]{or}}
\put(4,-4.5){\makebox(0,0)[lb]{\(\langle 4675\rangle\)}}


\end{picture}
}
\end{figure}

In January 2023,
we started a computationally extremely challenging
search for cyclic cubic fields \(K\) with HBC
and conductors \(10^5<c<3\cdot 10^6\).
Since the \(3\)-class tower group \(\mathrm{Gal}(\mathrm{F}_3^\infty(K)/K)\)
can only be identified by \textit{AQI of second order} (Table
\ref{tbl:AQI2HBC}),
we constructed the unramified cyclic cubic relative extensions \(E_{j,\ell}/E_j\)
of absolute degree \(27\)
and the class groups \(\mathrm{Cl}(E_{j,\ell})\),
\(1\le j\le 13\), \(1\le\ell\le n_j\),
forming the
\textit{Artin pattern \(\alpha^{(2)}\) of second order},
\begin{equation}
\label{eqn:AQI2NT}
\alpha^{(2)}:=\lbrack \mathrm{Cl}(E_j);(\mathrm{Cl}(E_{j,\ell}))_{\ell=1}^{n_j}\rbrack_{j=1}^{13}.
\end{equation}
Here,
\(1\le\ell\le 4\), \(n_j=4\), when \(\mathrm{Cl}_3(E_j)\simeq (2^2)\hat{=}(9,9)\), and
\(1\le\ell\le 13\), \(n_j=13\), when \(\mathrm{Cl}_3(E_j)\simeq (21^2)\hat{=}(9,3,3)\).
Since this was impossible even on bigger workstations,
due to the required RAM storage and CPU time,
we had to employ super computers with 128 cores
and 1TB RAM, which enabled highly parallel processes
with PARI/GP
\cite{PARI2023}
under the GRH.
\textit{Prototypes} with minimal conductors \(c\) (\textbf{boldface}) are visualized in Figure
\ref{fig:Realizations}.


\subsection{Galois structure of unramified cubic and nonic extensions}
\label{ss:Galois}

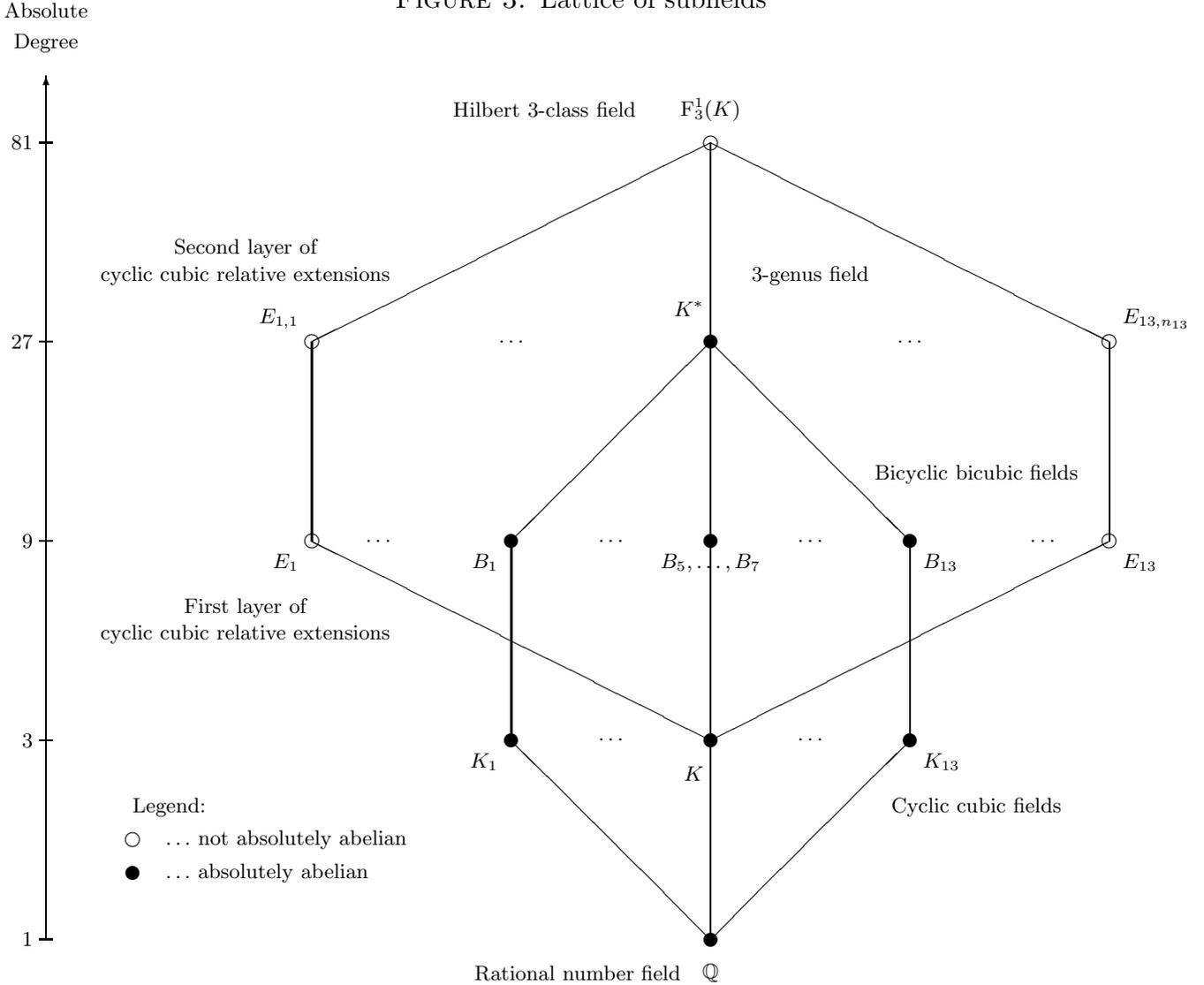
\begin{figure}[ht]
\caption{Lattice of subfields}
\label{fig:SubFieldLattice}

{\scriptsize

\setlength{\unitlength}{1.0cm}
\begin{picture}(15,14)(-9.5,-12)


\put(-10,2.0){\makebox(0,0)[cc]{Absolute}}
\put(-10,1.5){\makebox(0,0)[cc]{Degree}}

\put(-10,0){\vector(0,1){1}}

\put(-10,0){\line(0,-1){12}}
\multiput(-10.1,0)(0,-3){5}{\line(1,0){0.2}}

\put(-10.2,0){\makebox(0,0)[rc]{\(81\)}}
\put(-10.2,-3){\makebox(0,0)[rc]{\(27\)}}
\put(-10.2,-6){\makebox(0,0)[rc]{\(9\)}}
\put(-10.2,-9){\makebox(0,0)[rc]{\(3\)}}
\put(-10.2,-12){\makebox(0,0)[rc]{\(1\)}}


\put(-8.7,-10){\makebox(0,0)[lc]{Legend:}}
\put(-8.7,-10.5){\circle{0.2}}
\put(-8.2,-10.5){\makebox(0,0)[lc]{\(\ldots\) not absolutely abelian}}
\put(-8.7,-11){\circle*{0.2}}
\put(-8.2,-11){\makebox(0,0)[lc]{\(\ldots\) absolutely abelian}}


\put(-2.5,0.5){\makebox(0,0)[cc]{Hilbert \(3\)-class field}}
\put(0,0.5){\makebox(0,0)[cc]{\(\mathrm{F}_3^1(K)\)}}
\put(0,0){\circle{0.2}}
\put(0,0){\line(-2,-1){6}}
\put(0,0){\line(2,-1){6}}

\put(-7,-1.6){\makebox(0,0)[cc]{Second layer of}}
\put(-7,-2){\makebox(0,0)[cc]{cyclic cubic relative extensions}}
\multiput(-6,-3)(12,0){2}{\circle{0.2}}
\multiput(-3,-3)(6,0){2}{\makebox(0,0)[cc]{\(\cdots\)}}
\put(-6.2,-2.8){\makebox(0,0)[rb]{\(E_{1,1}\)}}
\put(6.2,-2.8){\makebox(0,0)[lb]{\(E_{13,n_{13}}\)}}
\put(-6,-6){\line(0,1){3}}
\put(6,-6){\line(0,1){3}}

\put(-7,-7){\makebox(0,0)[cc]{First layer of}}
\put(-7,-7.4){\makebox(0,0)[cc]{cyclic cubic relative extensions}}
\put(-6,-6){\circle{0.2}}
\put(-6.2,-6.2){\makebox(0,0)[rt]{\(E_1\)}}
\multiput(-5,-6)(10,0){2}{\makebox(0,0)[cc]{\(\cdots\)}}
\put(6,-6){\circle{0.2}}
\put(6.2,-6.2){\makebox(0,0)[lt]{\(E_{13}\)}}

\put(0,-9){\line(-2,1){6}}
\put(0,-9){\line(2,1){6}}
\put(0,-9){\circle*{0.2}}
\put(-0.1,-9.5){\makebox(0,0)[rc]{\(K\)}}


\put(1.5,-2){\makebox(0,0)[cc]{\(3\)-genus field}}
\put(-0.1,-2.5){\makebox(0,0)[rc]{\(K^\ast\)}}
\put(0,-3){\circle*{0.2}}
\put(0,-3){\line(-1,-1){3}}
\put(0,-3){\line(1,-1){3}}

\put(4,-5){\makebox(0,0)[cc]{Bicyclic bicubic fields}}
\put(-3,-6){\circle*{0.2}}
\put(-3.2,-6.2){\makebox(0,0)[rt]{\(B_1\)}}
\multiput(-1.5,-6)(3,0){2}{\makebox(0,0)[cc]{\(\cdots\)}}
\put(3,-6){\circle*{0.2}}
\put(3.2,-6.2){\makebox(0,0)[lt]{\(B_{13}\)}}
\put(0,-6){\circle*{0.2}}
\put(0,-6.2){\makebox(0,0)[ct]{\(B_5,\ldots,B_7\)}}
\put(-3,-9){\line(0,1){3}}
\put(3,-9){\line(0,1){3}}

\put(4,-10){\makebox(0,0)[cc]{Cyclic cubic fields}}
\put(-3,-9){\circle*{0.2}}
\put(-3.2,-9.2){\makebox(0,0)[rt]{\(K_1\)}}
\multiput(-1.5,-9)(3,0){2}{\makebox(0,0)[cc]{\(\cdots\)}}
\put(3,-9){\circle*{0.2}}
\put(3.2,-9.2){\makebox(0,0)[lt]{\(K_{13}\)}}

\put(0,-12){\line(-1,1){3}}
\put(0,-12){\line(0,1){12}}
\put(0,-12){\line(1,1){3}}
\put(-2,-12.5){\makebox(0,0)[cc]{Rational number field}}
\put(0,-12){\circle*{0.2}}
\put(0,-12.5){\makebox(0,0)[cc]{\(\mathbb{Q}\)}}

\end{picture}
}
\end{figure}

\noindent
The fundamental facts, on which the \textit{Galois structure}
of the lattice of intermediate fields
\(\mathbb{Q}<F<\mathrm{F}_3^1(K)\)
of the Hilbert \(3\)-class field \(\mathrm{F}_3^1(K)\)
of a cyclic cubic number field \(K\)
with \(3\)-class group \(\mathrm{Cl}_3(K)\simeq (3,3,3)\)
and conductor \(c=q_1q_2q_3\) with precisely three
prime (power) divisors \(q_i\equiv +1\,(\mathrm{mod}\,3)\),
or \(q_i=3^2\), is based, can be summarized as follows (see Figure
\ref{fig:SubFieldLattice}):
\begin{itemize}
\item
The \textit{absolute \(3\)-genus field} \(K^\ast=(K/\mathbb{Q})^\ast\) of \(K\)
is the maximal unramified \(3\)-exten\-sion of \(K\)
which is abelian over the rational field \(\mathbb{Q}\).
In the situation with conductor \(c=q_1q_2q_3\),
its absolute Galois group is
\(\mathrm{Gal}(K^\ast/\mathbb{Q})\simeq (3,3,3)\),
whence it possesses
\(13\) cyclic cubic subfields \(K_1,\ldots,K_{13}\),
one of them \(K\), and
\(13\) bicyclic bicubic subfields \(B_1,\ldots,B_{13}\).
The former consist of
three singlets, three doublets, and a quartet
with partial conductors \(q_1,q_2,q_3\),
resp. \(q_1q_2,q_2q_3,q_3q_1\),
resp. \(c=q_1q_2q_3\).
The exact constitution of the latter was analyzed by Ayadi
\cite{Ay1995}:
Three of them are sub-genus fields
\(B_{11}=k_{q_1q_2}^\ast,B_{12}=k_{q_2q_3}^\ast,B_{13}=k_{q_3q_1}^\ast\)
with conductors \(q_1q_2,q_2q_3,q_3q_1\).
Among the remaining ten, four contain \(K\), namely, in Ayadi's notation
\cite[p. 42]{Ay1995},\\
\(B_1=Kk_{q_1q_2}k_{q_1q_3}k_{q_2q_3}\),
\(B_5=KK_3k_{q_1}\tilde{k}_{q_2q_3}\),
\(B_6=KK_4k_{q_2}\tilde{k}_{q_1q_3}\),
\(B_7=KK_2k_{q_3}\tilde{k}_{q_1q_2}\).
\end{itemize}


\begin{itemize}
\item
The \textit{Hilbert \(3\)-class field} \(\mathrm{F}_3^1(K)\) of \(K\)
is the maximal abelian unramified \(3\)-extension of \(K\).
By Artin's reciprocity law, its relative Galois group
\(\mathrm{Gal}(\mathrm{F}_3^1(K)/K)\) is isomorphic to
the \(3\)-class group \(\mathrm{Cl}_3(K)\simeq  (3,3,3)\) of \(K\).
Among the \(13\) cyclic cubic relative extensions
\(K<E_1,\ldots,E_{13}<\mathrm{F}_3^1(K)\),
only four are abelian, namely the bicyclic bicubic fields
\(B_1,B_5,B_6,B_7\)
with absolute Galois group \(\langle 9,2\rangle\).
The remaining nine extensions are non-Galois with
Galois group \(\langle 27,3\rangle\)
of the splitting field,
arranged in three triplets of three isomorphic fields each.
\end{itemize}


\noindent
An unsophisticated way to determine the Artin pattern \(\alpha^{(2)}\) of second order
of a cyclic cubic field \(K\)
would be to construct the entire collection of
the following extensions \(E_{j,\ell}\).

Recall that we have \textit{three possible scenarios}, according to Table
\ref{tbl:HBC}:
\begin{enumerate}
\item
the rank distribution \(3^1,2^{12}\),
equivalently the Taussky types \(\mathrm{A}^1,\mathrm{B}^{12}\),
with \(1\cdot 13+12\cdot 4=13+48=61\) unramified nonic
but not necessarily Galois extensions \(E_{j,\ell}/K\);
\item
the rank distribution \(3^4,2^9\),
equivalently the Taussky types \(\mathrm{A}^4,\mathrm{B}^9\),
with \(4\cdot 13+9\cdot 4=52+36=88\) unramified nonic
but not necessarily Galois extensions \(E_{j,\ell}/K\);
\item
the rank distribution \(3^7,2^6\),
equivalently the Taussky types \(\mathrm{A}^7,\mathrm{B}^6\),
with \(7\cdot 13+6\cdot 4=91+24=115\) unramified nonic
but not necessarily Galois extensions \(E_{j,\ell}/K\).
\end{enumerate}

\noindent
We avoid the computation of \(61\), respectively \(88\), respectively \(115\),
extensions \(E_{j,\ell}\) and their \(3\)-class groups
\(\mathrm{Cl}_3(E_{j,\ell})\)
by using \textit{isomorphisms to representatives} \(R_i\simeq E_{j,\ell}\).

Firstly, we only need \textit{seven} extensions (Table
\ref{tbl:FirstLayer}),
the four abelian \(B_1,B_5,B_6,B_7\) and
three non-Galois \(E_j\), one of each triplet of three isomorphic fields,
in the \textit{first layer} of unramified cyclic cubic relative extensions
\(E_1,\ldots,E_{13}\) of \(K\),
which are of absolute degree \(9\).
Their \(3\)-class groups \(\lbrack\mathrm{Cl}_3(E_j)\rbrack_{j=1}^{13}\)
constitute the Artin pattern \(\alpha(K)=\alpha^{(1)}(K)\) of \textit{first order} of \(K\). \\
(In the column \(\#\), the symbol \(n/m\) denotes \(n\) conjugacy classes with \(m\) members each.)


\renewcommand{\arraystretch}{1.1}

\begin{table}[ht]
\caption{Isomorphisms and representatives among extensions of degree \(9\)}
\label{tbl:FirstLayer}
\begin{center}
{\scriptsize
\begin{tabular}{|c|c||c||c|c|c||c|c|c||c|c|}
\hline
 Sc. & \(\rho\)       &    & \(\#\) & Rep.            & Abelian                          & \(\#\)  & Rep.            & Non-Galois                        &      & Census         \\
\hline
 (1) & \(3^1,2^{12}\) &  1 & \(1\)  & \(B_1\)         & \(21^2\), \(\langle 9,2\rangle\) &         &                 &                                   & Rep. & \(1+3+3=7\)    \\
     &                & 12 & \(3\)  & \(B_5,B_6,B_7\) & \(2^2\),  \(\langle 9,2\rangle\) & \(3/3\) & \(E_1,E_2,E_3\) & \(2^2\),  \(\langle 27,3\rangle\) & Tot. & \(1+3+9=13\)   \\
\hline
 (2) & \(3^4,2^9\)    &  4 & \(1\)  & \(B_1\)         & \(21^2\), \(\langle 9,2\rangle\) & \(1/3\) & \(E_1\)         & \(21^2\), \(\langle 27,3\rangle\) & Rep. & \(1+1+3+2=7\)  \\
     &                &  9 & \(3\)  & \(B_5,B_6,B_7\) & \(2^2\),  \(\langle 9,2\rangle\) & \(2/3\) & \(E_2,E_3\)     & \(2^2\),  \(\langle 27,3\rangle\) & Tot. & \(1+3+3+6=13\) \\
\hline
 (3) & \(3^7,2^6\)    &  7 & \(1\)  & \(B_1\)         & \(21^2\), \(\langle 9,2\rangle\) & \(2/3\) & \(E_1,E_2\)     & \(21^2\), \(\langle 27,3\rangle\) & Rep. & \(1+2+3+1=7\)  \\
     &                &  6 & \(3\)  & \(B_5,B_6,B_7\) & \(2^2\),  \(\langle 9,2\rangle\) & \(1/3\) & \(E_3\)         & \(2^2\),  \(\langle 27,3\rangle\) & Tot. & \(1+6+3+3=13\) \\
\hline
\end{tabular}
}
\end{center}
\end{table}


\noindent
Secondly, among the fields of absolute degree \(27\)
in the \textit{second layer}, we can restrict our class group computations to
\begin{enumerate}
\item
\textit{eight} for the first scenario with rank distribution \((3^1,2^{12})\),
\item
\textit{eleven} for the second scenario with rank distribution \((3^4,2^9)\), and
\item
\textit{fourteen} for the third scenario with rank distribution \((3^7,2^6)\) (Table
\ref{tbl:SecondLayer}).
\end{enumerate}
For each scenario, only the \(3\)-genus field \(K^\ast=:R_1\) is abelian
with absolute Galois group \(\langle 27,5\rangle\),
the other fields may be Galois with group \(\langle 27,3\rangle\)
or non-Galois with various groups of the splitting field:
\(\langle 729,411\rangle\) for \(\rho=(3^4,2^9)\) or \(\rho=(3^7,2^6)\), and
\(\langle 81,12\rangle\), \(\langle 243,58\rangle\) for all scenarios.
Their \(3\)-class groups 
\(\lbrack (\mathrm{Cl}_3(E_{j,\ell}))_{\ell=1}^{n_j}\rbrack_{j=1}^{13}\)
constitute the Artin pattern \(\alpha^{(2)}(K)\) of \textit{second order} of \(K\).


\renewcommand{\arraystretch}{1.1}

\begin{table}[ht]
\caption{Isomorphisms and representatives among extensions of degree \(27\)}
\label{tbl:SecondLayer}
\begin{center}
{\scriptsize
\begin{tabular}{|c|c||c||c|c|c||c|c|c||c|c|c|}
\hline
 Sc. & \(\rho\)       &       & \(\#\) & Rep.    & Abelian                 & \(\#\) & Rep.            & Galois                  & \(\#\)  &  Rep.                 & Non-Galois                \\
\hline
 (1) & \(3^1,2^{12}\) & \(1\) & \(1\)  & \(R_1\) & \(\langle 27,5\rangle\) & \(3\)  & \(R_2,R_3,R_4\) & \(\langle 27,3\rangle\) & \(1/9\) & \(R_8\)               & \(\langle 243,58\rangle\) \\
     &                & \(3\) & \(3\)  & \(R_1\) & \(\langle 27,5\rangle\) &        &                 &                         & \(3/3\) & \(R_5,R_6,R_7\)       & \(\langle 81,12\rangle\)  \\
     &                & \(9\) &        &         &                         & \(9\)  & \(R_2,R_3,R_4\) & \(\langle 27,3\rangle\) & \(9/3\) & \(R_5,R_6,R_7\)       & \(\langle 81,12\rangle\)  \\
     & Subtotal       &       & \(4\)  &         &                         & \(12\) &                 &                         & \(45\)  & Total                 & \(4+12+45=61\)            \\
\hline
 (2) & \(3^4,2^9\)    & \(1\) & \(1\)  & \(R_1\) & \(\langle 27,5\rangle\) & \(3\)  & \(R_2,R_3,R_4\) & \(\langle 27,3\rangle\) & \(1/9\) & \(R_8\)               & \(\langle 243,58\rangle\) \\
     &                & \(3\) &        &         &                         & \(3\)  & \(R_2\)         & \(\langle 27,3\rangle\) & \(3/3\) & \(R_5,R_6,R_7\)       & \(\langle 81,12\rangle\)  \\
     &                &       &        &         &                         &        &                 &                         & \(3/9\) & \(R_9,R_{10},R_{11}\) & \(\langle 729,411\rangle\)\\
     &                & \(3\) & \(3\)  & \(R_1\) & \(\langle 27,5\rangle\) &        &                 &                         & \(3/3\) & \(R_5,R_6,R_7\)       & \(\langle 81,12\rangle\)  \\
     &                & \(6\) &        &         &                         & \(6\)  & \(R_3,R_4\)     & \(\langle 27,3\rangle\) & \(6/3\) & \(R_5,R_6,R_7\)       & \(\langle 81,12\rangle\)  \\
     & Subtotal       &       & \(4\)  &         &                         & \(12\) &                 &                         & \(72\)  & Total                 & \(4+12+72=88\)            \\
\hline
 (3) & \(3^7,2^6\)    & \(1\) & \(1\)  & \(R_1\) & \(\langle 27,5\rangle\) & \(3\)  & \(R_2,R_3,R_4\) & \(\langle 27,3\rangle\) & \(1/9\) & \(R_8\)               & \(\langle 243,58\rangle\) \\
     &                & \(3\) &        &         &                         & \(3\)  & \(R_2\)         & \(\langle 27,3\rangle\) & \(3/3\) & \(R_5,R_6,R_7\)       & \(\langle 81,12\rangle\)  \\
     &                &       &        &         &                         &        &                 &                         & \(3/9\) & \(R_9,R_{10},R_{11}\) & \(\langle 729,411\rangle\)\\
     &                & \(3\) &        &         &                         & \(3\)  & \(R_3\)         & \(\langle 27,3\rangle\) & \(3/3\) & \(R_5,R_6,R_7\)       & \(\langle 81,12\rangle\)  \\
     &                &       &        &         &                         &        &                 &                         & \(3/9\) &\(R_{12},R_{13},R_{14}\)&\(\langle 729,411\rangle\)\\
     &                & \(3\) & \(3\)  & \(R_1\) & \(\langle 27,5\rangle\) &        &                 &                         & \(3/3\) & \(R_5,R_6,R_7\)       & \(\langle 81,12\rangle\)  \\
     &                & \(3\) &        &         &                         & \(3\)  & \(R_4\)         & \(\langle 27,3\rangle\) & \(3/3\) & \(R_5,R_6,R_7\)       & \(\langle 81,12\rangle\)  \\
     & Subtotal       &       & \(4\)  &         &                         & \(12\) &                 &                         & \(99\)  & Total                 & \(4+12+99=115\)           \\
\hline
\end{tabular}
}
\end{center}
\end{table}


\section{Computations}
\label{s:Computations}

The computations were performed using the PARI/GP\cite{PARI2023} computer algebra system.
The two most important steps are the computations of class groups and unit groups
(performed by the GP function \texttt{bnfinit})
and the computations of class fields (performed by the GP function \texttt{bnrclassfield}).
A call to \texttt{bnrclassfield} for a field uses Kummer theory and requires
calls to \texttt{bnfinit} for both the field and its extension by the third roots of unity (in our case).
So the actual computation consists of these steps, starting from a suitable cyclic cubic field:
\begin{enumerate}
\item compute the class group and unit group of the cubic field
\item compute the class group and unit group of the cubic field extended with the $3$-root of unity.
\item compute class fields of degree $9$
\item for each class field,
\begin{enumerate}
\item compute its class group and unit group
\item compute the class group and unit group of the field extended with the $3$-root of unity.
\item compute class fields of degree $27$
\item compute the class groups of each class fields.
\end{enumerate}
\end{enumerate}
As a result, for each cyclic cubic field, we need to call \texttt{bnfinit} on
one field of degree $3$, one field of degree $6$, $7$ fields of degree $9$, $7$
fields of degree $18$, and either $8$, $11$ or $14$ non-isomorphic fields of
degree $27$.

The function \texttt{bnfinit} is an implementation of Buchmann
subexponential algorithm for class groups and unit groups~\cite{Bu1990}
by Cohen-Diaz-Olivier~\cite{CDO1997},~\cite{Co2000}. It is based on searching 
relations between ideals in a set of prime ideals that generates the class group,
and is correct under the assumption of the Riemann hypothesis for all Hecke
$L$-functions attached to non-trivial characters of the ideal class group~\cite{BeFm2015}.

This computation was done with a specially-tuned, parallel version of this function.
The set of primes is chosen by applying Greni\'e-Molteni~\cite{GnMo2016} improvement to
Belabas-Diaz-Friedman \cite{BDF2008} criterion to find a small set of primes generating
the class group.
It proceeds by searching in parallel for smooth elements in the ideals 
obtained by applying LLL-reduction to the ideals $\mathfrak{p}^6\mathfrak{q}$ for all pairs of ideals
$(\mathfrak{p},\mathfrak{q})$ in the set. The tuning parameters decide how far the program will search
in each such ideal. It was regularly increased to account for the increase in the fields discriminant
over the course of the computation.
While the units are not required for the fields of degree~$27$, we still used the compact units
representation of units because precision increases would be parallelised and were less expensive
than with the logarithmic embedding representation.
The program was run on a 128-core CPU with 1TB of RAM over the course of several months, using the
internal POSIX threads parallel engine of PARI/GP.


\section{Technical details}
\label{s:Details}

\noindent
The following details are not required to understand
our main results on cyclic cubic fields \(K\) with HBC.
They are, however,
mandatory for the correct selection of
representatives \(R_i\) in isomorphism classes
among the unramified cyclic cubic relative extensions
\(E_j/K\), \(1\le j\le 13\),
and among the unramified nonic but not necessarily Galois extensions
\(E_{j,\ell}/K\), \(1\le j\le 13\), \(1\le\ell\le n_j\), \(n_j\in\lbrace 4,13,40\rbrace\),
of absolute degree \(27\).

\begin{theorem}
\label{thm:CycCub}
Let \(K\) be a cyclic cubic number field
with conductor \(c=q_1q_2q_3\) divisible by exactly three distinct prime(power)s,
\(q_i\equiv +1\,(\mathrm{mod}\,3)\), or \(q_i=3^2\).
Then
\begin{enumerate}
\item
\(K\) is member of a quartet \((K_1,\ldots,K_4)\) of four cyclic cubic fields
sharing the common conductor \(c\) and the common discriminant \(d=c^2\).
\item
The absolute genus field \(K^\ast=(K/\mathbb{Q})^\ast\) of \(K\)
is unramified over \(K\) and abelian over \(\mathbb{Q}\). More precisely,
its Galois group \(\mathrm{Gal}(K^\ast/\mathbb{Q})\simeq(\mathbb{Z}/3\mathbb{Z})^3\)
is elementary tricyclic. 
Its absolute degree is \(\lbrack K^\ast:\mathbb{Q}\rbrack=27\) and
the relative degree is \(\lbrack K^\ast:K\rbrack=9\).
\item
\(K^\ast\) contains \(13\) cyclic cubic subfields,
three \(k_{q_1},k_{q_2},k_{q_3}\) with prime(power) conductors,
six (in three doublets)
\(k_{q_1q_2},\tilde{k}_{q_1q_2},k_{q_1q_3},\tilde{k}_{q_1q_3},k_{q_2q_3},\tilde{k}_{q_2q_3}\)
with conductors divisible by two prime(power)s,
and the four members \(K_1,\ldots,K_4\) of the abovementioned quartet with conductor \(c\).
\item
The composita \(L:=k_{q_1q_2}k_{q_1q_3}k_{q_2q_3}\) and
\(\tilde{L}:=\tilde{k}_{q_1q_2}\tilde{k}_{q_1q_3}\tilde{k}_{q_2q_3}\)
satisfy the following \textbf{skew balance of degrees}:
\(\lbrack L:\mathbb{Q}\rbrack\cdot\lbrack\tilde{L}:\mathbb{Q}\rbrack=243\), with
\begin{equation}
\label{eqn:Skew}
\lbrack L:\mathbb{Q}\rbrack=9 \Longleftrightarrow \lbrack\tilde{L}:\mathbb{Q}\rbrack=27,
\end{equation}
or vice versa.
\end{enumerate}
\end{theorem}

\begin{proof}
See
\cite[\S\ 4.1, p. 40, Proof of Prop. 4.6, p. 49, Prop. 4.1, p. 40]{Ay1995}.
\end{proof}

\begin{definition}
\label{dfn:Normalization}
The selection of cyclic cubic subfields \(k_{q_1q_2},k_{q_1q_3},k_{q_2q_3}\) with
conductors \(q_1q_2,q_1q_3,q_2q_3\) within the absolute genus field \(K^\ast\)
of an assigned cyclic cubic field \(K\) with conductor \(c=q_1q_2q_3\)
is called \textit{normalized}, if the absolute degree of their compositum
\(L=k_{q_1q_2}k_{q_1q_3}k_{q_2q_3}\) is \(\lbrack L:\mathbb{Q}\rbrack=9\).
In this article, we always assume this normalization.
\end{definition}

\begin{theorem}
\label{thm:BicycBicub}
Under the assumptions of Theorem
\ref{thm:CycCub}
and the mandatory normalization of \(k_{q_1q_2},k_{q_1q_3},k_{q_2q_3}\)
according to Definition
\ref{dfn:Normalization},
the remaining \(13\) bicyclic bicubic subfields \(B_j\), \(1\le j\le 13\),
of the absolute genus field \(K^\ast\) of \(K\) are given as composita by
\begin{equation}
\label{eqn:BicycBicub}
\begin{aligned}
B_1 &:= k_{q_1q_2}k_{q_1q_3} = K_1k_{q_1q_2}k_{q_1q_3}k_{q_2q_3} \\
B_2 &:= \tilde{k}_{q_1q_3}\tilde{k}_{q_2q_3} = K_2k_{q_1q_2}\tilde{k}_{q_1q_3}\tilde{k}_{q_2q_3} \\
B_3 &:= \tilde{k}_{q_1q_2}\tilde{k}_{q_1q_3} = K_3\tilde{k}_{q_1q_2}\tilde{k}_{q_1q_3}k_{q_2q_3} \\
B_4 &:= \tilde{k}_{q_1q_2}\tilde{k}_{q_2q_3} = K_4\tilde{k}_{q_1q_2}k_{q_1q_3}\tilde{k}_{q_2q_3} \\
B_5 &:= K_1K_3 = K_1K_3k_{q_1}\tilde{k}_{q_2q_3} \\
B_6 &:= K_1K_4 = K_1K_4k_{q_2}\tilde{k}_{q_1q_3} \\
B_7 &:= K_1K_2 = K_1K_2k_{q_3}\tilde{k}_{q_1q_2} \\
B_8 &:= K_2K_4 = K_2K_4k_{q_1}k_{q_2q_3} \\
B_9 &:= K_2K_3 = K_2K_3k_{q_2}k_{q_1q_3} \\
B_{10} &:= K_3K_4 = K_3K_4k_{q_3}k_{q_1q_2} \\
B_{11} &:= k_{q_1q_2}\tilde{k}_{q_1q_2} = k_{q_1}k_{q_2}k_{q_1q_2}\tilde{k}_{q_1q_2} \\
B_{12} &:= k_{q_1q_3}\tilde{k}_{q_1q_3} = k_{q_1}k_{q_3}k_{q_1q_3}\tilde{k}_{q_1q_3} \\
B_{13} &:= k_{q_2q_3}\tilde{k}_{q_2q_3} = k_{q_2}k_{q_3}k_{q_2q_3}\tilde{k}_{q_2q_3}
\end{aligned}
\end{equation}
The shape with two components suffices for the construction,
but the shape with four components ostensively illuminates
all cyclic cubic subfields of each bicyclic bicubic field \(B_j\).
\end{theorem}

\begin{proof}
See Ayadi's Thesis
\cite[Lem. 4.1, p. 42, and Fig. 10, p. 41]{Ay1995}.
\end{proof}

\begin{corollary}
\label{cor:BicycBicub}
For each member of the quartet \((K_1,\ldots,K_4)\) of cyclic cubic fields
with conductor \(c=q_1q_2q_3\),
the rank \(\varrho_i\) of the \(3\)-class group \(\mathrm{Cl}_3(K_i)\) is bounded by
\(2\le\varrho_i\le 4\),
and four unramified cyclic cubic relative extensions of \(K_i\) are given in the following way:
\begin{equation}
\label{eqn:Unramified}
\begin{aligned}
(B_1,B_5,B_6,B_7) &\text{ for } K_1, \\
(B_2,B_7,B_8,B_9) &\text{ for } K_2, \\
(B_3,B_5,B_9,B_{10}) &\text{ for } K_3, \\
(B_4,B_6,B_8,B_{10}) &\text{ for } K_4. \\
\end{aligned}
\end{equation}
If the rank of the \(3\)-class group \(\mathrm{Cl}_3(K_i)\) of \(K_i\) is \(\varrho_i=2\),
then the set of unramified extensions given in Equation
\eqref{eqn:Unramified}
is complete and consists entirely of absolutely abelian extensions.
\end{corollary}

\begin{proof}
This is an immediate consequence of the constitution of the \(B_j\) in Theorem
\ref{thm:BicycBicub}.
\end{proof}


\section{Conclusion}
\label{s:Conclusion}

\noindent
Generally, our investigation of the \(3\)-class field tower of cyclic cubic fields \(K\)
with \textit{elementary tricyclic} \(3\)-class group \(\mathrm{Cl}_3(K)\simeq (3,3,3)\)
is a striking novelty
\cite{Ma2022}.
Similar attempts with imaginary quadratic fields of type \((3,3,3)\),
where all capitulation kernels are of order \(\#L=3\) (lines),
successfully yielded the Artin pattern \(\mathrm{AP}=(\alpha,\varkappa)\)
by means of arithmetic computations
\cite[\S\ 7.2, Tbl. 2--4, pp. 308--311]{Ma2015b}
but were doomed to failure group theoretically,
since the order of relevant groups is at least \(3^{15}\)
and the complexity of descendant trees became unmanageable
\cite[\S\ 7.4, p. 312]{Ma2015b},
\cite[\S\ 10, p. 54]{Ma2015c},
\cite[Thm. 8.2, p. 174]{Ma2017},
\cite[\S\ 8, pp. 98--99]{Ma2016c},
\cite[\S\ 2, Example, p. 6]{Ma2016d}.
Therefore, we were delighted that cyclic cubic fields of type \((3,3,3)\)
impose much less severe requirements on the second \(3\)-class group
\(\mathfrak{M}=\mathrm{Gal}(\mathrm{F}_3^2(K)/K)\),
since capitulation kernels of order \(\#P=9\) (planes)
and even \(\#O=27\) (full space) are admissible.
However, in the present work our attention is devoted to cyclic cubic fields
with \textit{harmonically balanced capitulation} (HBC),
where all transfer kernels are of order \(\#L=3\) (lines),
but relevant groups set in at order \(3^6\) already (Tbl.
\ref{tbl:HBC}).
Our foremost target was the realization of \textit{closed Andozhskii-Tsvetkov groups}
(AT-groups) \(G\)
with coinciding generator- and relation-rank \(d_2(G)=d_1(G)\).
The main result was:

\begin{theorem}
\label{thm:Cat1Gph2}
Let \(K\) be a cyclic cubic number field with conductor \(c=q_1q_2q_3\)
divisible by precisely three distinct prime(power)s,
\(q_i\equiv +1\,(\mathrm{mod}\,3)\), or \(q_i=3^2\),
such that only two cubic residue symbols
\(\left(\frac{q_1}{q_2}\right)_3=1\) and \(\left(\frac{q_1}{q_3}\right)_3=1\)
are trivial, that is, \(K\) belongs to \textbf{graph} \(\mathbf{2}\) of \textbf{category I},
\(q_2\leftarrow q_1\rightarrow q_3\),
in the sense of G. Gras and Ayadi. Then:
\begin{enumerate}
\item
The \(3\)-class group \(\mathrm{Cl}_3(K)\) has either rank \(\varrho=3\),
for a single component,
or it is elementary bicyclic
\(\mathrm{Cl}_3(K)\simeq(\mathbb{Z}/3\mathbb{Z})^2\), with \(\varrho=2\),
for three components of the quartet \((K_1,\ldots,K_4)\).
The former is \(K_1\) if
\(q_2\) splits in \(k_{q_1q_3}\) and
\(q_3\) splits in \(k_{q_1q_2}\).
\item
If \(\varrho=2\), then \(q_1\) is the unique minimal norm of 
a non-trivial primitive ambiguous principal ideal of \(K\),
called \textbf{Parry invariant}
\cite[pp. 499--501]{Pa1990}
of \(K\) by Ayadi.
\item
If
\(q_2\) splits in \(k_{q_1q_3}\),
\(q_3\) splits in \(k_{q_1q_2}\),
and \(K_1\) possesses an elementary tricyclic \(3\)-class group
\(\mathrm{Cl}_3(K)\simeq(\mathbb{Z}/3\mathbb{Z})^3\)
with \textbf{HBC}, then the Parry invariant of \(K_1\) is also \(q_1\)
and the remaining three fields \(K_2,K_3,K_4\) with \(\varrho=2\)
share the common capitulation type
\(\varkappa(K_i)\sim (1243)\) with
two fixed points \(1,2\) and a transposition \((43)\),
called type \(\mathrm{G}.16\),
and their \(3\)-class field tower has two stages with automorphism group
\(\langle 729,52\rangle\).
\end{enumerate}
\end{theorem}

\begin{proof}
Denote by \(G:=\mathrm{Gal}(K/\mathbb{Q})=\langle\sigma\rangle\)
the cyclic Galois group of \(K\).
Among the prime ideals of \(K\), let
\(\mathfrak{P}=\mathfrak{P}^\sigma\) divide \(q_1\),
\(\mathfrak{Q}=\mathfrak{Q}^\sigma\) divide \(q_2\), and
\(\mathfrak{R}=\mathfrak{R}^\sigma\) divide \(q_3\).
If the rank \(\varrho\) of the \(3\)-class group \(\mathrm{Cl}_3(K)\) is \(\varrho=2\),
then \(\mathrm{Cl}_3(K)\simeq(\mathbb{Z}/3\mathbb{Z})^2\)
\cite[Prop. 4.3, p. 43]{Ay1995}.
If \(\varrho=2\), then \(\mathfrak{P}\) generates the group
\(\mathcal{P}_K^G/\mathcal{P}_{\mathbb{Q}}\)
of primitive ambiguous principal ideals of \(K\)
\cite[Rem. 4.2, p. 50]{Ay1995},
whereas \(\mathfrak{Q}\) and \(\mathfrak{R}\) are not principal
\cite[Rem. 4.8, p. 65]{Ay1995},
and their ideal classes
\(\lbrack\mathfrak{Q}\rbrack=\mathfrak{Q}\cdot\mathcal{P}_K\) and
\(\lbrack\mathfrak{R}\rbrack=\mathfrak{R}\cdot\mathcal{P}_K\)
generate \(\mathrm{Cl}_3(K)=
\langle\lbrack\mathfrak{Q}\rbrack,\lbrack\mathfrak{R}\rbrack\rangle\).
If \(q_2\) splits in \(k_{q_1q_3}\) and
\(q_3\) splits in \(k_{q_1q_2}\),
then \(K_1\) is the field with \(\varrho=3\)
\cite[Prop. 4.4, pp. 43--44]{Ay1995},
and \(K_2,K_3,K_4\) have elementary bicyclic \(3\)-class groups.
According to
\cite[Tbl., p. 66]{Ay1995},
the kernels of the transfers
\(T_{ij}:\,\mathrm{Cl}_3(K_i)\to\mathrm{Cl}_3(B_j)\)
from \(K_i\), \(2\le i\le 4\), to its four unramified cyclic cubic extensions \(B_j\),
given in Corollary
\ref{cor:BicycBicub},
are as follows:
\begin{equation}
\label{eqn:Capitulation}
\begin{aligned}
\ker(T_{22})&=\langle\lbrack\mathfrak{Q}\mathfrak{R}^2\rbrack\rangle,\
\ker(T_{27})=\langle\lbrack\mathfrak{R}\rbrack\rangle,\
\ker(T_{28})=\langle\lbrack\mathfrak{Q}\mathfrak{R}\rbrack\rangle,\
\ker(T_{29})=\langle\lbrack\mathfrak{Q}\rbrack\rangle; \\
\ker(T_{33})&=\langle\lbrack\mathfrak{Q}\mathfrak{R}\rbrack\rangle,\
\ker(T_{35})=\langle\lbrack\mathfrak{Q}\mathfrak{R}^2\rbrack\rangle,\
\ker(T_{39})=\langle\lbrack\mathfrak{Q}\rbrack\rangle,\
\ker(T_{3,10})=\langle\lbrack\mathfrak{R}\rbrack\rangle; \\
\ker(T_{44})&=\langle\lbrack\mathfrak{Q}\mathfrak{R}^2\rbrack\rangle,\
\ker(T_{46})=\langle\lbrack\mathfrak{Q}\rbrack\rangle,\
\ker(T_{48})=\langle\lbrack\mathfrak{Q}\mathfrak{R}\rbrack\rangle,\
\ker(T_{4,10})=\langle\lbrack\mathfrak{R}\rbrack\rangle.
\end{aligned}
\end{equation}
For each row \(2\le i\le 4\), the transfer kernels form a permutation
of the four cyclic subgroups of order \(3\) of \(\mathrm{Cl}_3(K)=
\langle\lbrack\mathfrak{Q}\rbrack,\lbrack\mathfrak{R}\rbrack\rangle\),
more precisely, each row has two fixed points,
where the norm class group \(N_{B_j/K_i}(\mathrm{Cl}_3(K_i))\)
coincides with the transfer kernel \(\ker(T_{ij})\),
and a transposition,
where the norm class group and the transfer kernel are twisted.
This characterizes type \(\mathrm{G}.16\) unambiguously,
and the Galois group \(\mathrm{Gal}(F_3^\infty(K_i)/K_i)\)
of the \(3\)-class field tower is
the metabelian \(3\)-group \(S=\langle 729,52\rangle\)
with coclass \(\mathrm{cc}(S)=2\) and relation rank \(d_2(S)=4\)
as required for \(\varrho=2\),
for \(2\le i\le 4\).
\end{proof}


\section{Acknowledgements}
\label{s:Acknowledgements}

\noindent
Experiments presented in this paper were carried out using the
PlaFRIM
(Plateforme F\'ed\'erative pour la Recherche en Informatique et Math\'ematique)
experimental testbed,
supported by Inria,
CNRS (LABRI and IMB),
Universit\'e de Bordeaux,
Bordeaux INP and
Conseil R\'egional d'Aquitaine
(see
\texttt{https://www.plafrim.fr}).
All the computations were carried out on a system running
two 64-core AMD Zen3 EPYC 7763 CPUs
at 2.45 GHz
with 1TB of RAM.

The second author acknowledges that his research was supported by
the Austrian Science Fund (FWF): projects J0497-PHY, P26008-N25,
and by the Research Executive Agency of the European Union (EUREA):
project Horizon Europe 2021--2027.



\end{document}